\documentclass[11pt,reqno]{amsart}

\usepackage{amsthm}
\usepackage{amssymb}
\usepackage{latexsym}
\usepackage{multicol}
\usepackage{verbatim,enumerate}
\usepackage{hyperref}
\usepackage{amsmath, amscd}
\usepackage[all,cmtip]{xy}
\usepackage{soul}

\advance\textwidth by 1.2in \advance\oddsidemargin by -.6in \advance\evensidemargin by -.6in
\usepackage{tikz}
\usetikzlibrary{decorations.markings}
\tikzstyle{vertex}=[ circle, fill, draw, inner sep=0pt, minimum size=4pt,]
\tikzstyle{edge}= [thick]

\newtheorem*{lem}{Lemma}
\newtheorem*{prop}{Proposition}

\theoremstyle{definition} 
\theoremstyle{definition}
\newtheorem{thm}{Theorem}
\newtheorem*{thm*}{Theorem}

\newtheorem*{rem}{Remark}

\newenvironment{pf}{\proof}{\endproof}
\newcounter{cnt}
\newenvironment{enumerit}{\begin{list}{{\hfill\rm(\roman{cnt})\hfill}}{%
\settowidth{\labelwidth}{{\rm(iv)}}\leftmargin=\labelwidth%
\advance\leftmargin by \labelsep\rightmargin=0pt\usecounter{cnt}}}{\end{list}} \makeatletter
\def\mydggeometry{\makeatletter\dg@YGRID=1\dg@XGRID=20\unitlength=0.003pt\makeatother}
\makeatother \theoremstyle{remark}

\numberwithin{equation}{section}

 \DeclareMathOperator{\Ht}{ht}

\begin{document}

\newcommand{\thmref}[1]{Theorem~\ref{#1}}
\newcommand{\secref}[1]{Section~\ref{#1}}
\newcommand{\lemref}[1]{Lemma~\ref{#1}}
\newcommand{\propref}[1]{Proposition~\ref{#1}}
\newcommand{\corref}[1]{Corollary~\ref{#1}}
\newcommand{\remref}[1]{Remark~\ref{#1}}
\newcommand{\defref}[1]{Definition~\ref{#1}}
\newcommand{\er}[1]{(\ref{#1})}
\newcommand{\id}{\operatorname{id}}
\newcommand{\ord}{\operatorname{\emph{ord}}}
\newcommand{\sgn}{\operatorname{sgn}}
\newcommand{\wt}{\operatorname{wt}}
\newcommand{\tensor}{\otimes}
\newcommand{\from}{\leftarrow}
\newcommand{\nc}{\newcommand}
\newcommand{\rnc}{\renewcommand}
\newcommand{\dist}{\operatorname{dist}}
\newcommand{\qbinom}[2]{\genfrac[]{0pt}0{#1}{#2}}
\nc{\cal}{\mathcal} \nc{\goth}{\mathfrak} \rnc{\bold}{\mathbf}
\renewcommand{\frak}{\mathfrak}
\newcommand{\supp}{\operatorname{supp}}
\newcommand{\Irr}{\operatorname{Irr}}
\newcommand{\psym}{\mathcal{P}^+_{K,n}}
\newcommand{\psyml}{\mathcal{P}^+_{K,\lambda}}
\newcommand{\psymt}{\mathcal{P}^+_{2,\lambda}}
\renewcommand{\Bbb}{\mathbb}
\nc\bomega{{\mbox{\boldmath $\omega$}}} \nc\bpsi{{\mbox{\boldmath $\Psi$}}}
 \nc\balpha{{\mbox{\boldmath $\alpha$}}}
 \nc\bbeta{{\mbox{\boldmath $\beta$}}}
 \nc\bpi{{\mbox{\boldmath $\pi$}}}
  \nc\bpis{{\mbox{\boldmath \scriptsize$\pi$}}}
 \nc\bullets{{\mbox{\scriptsize $\bullet$}}}
 
  \nc\bvarpi{{\mbox{\boldmath $\varpi$}}}

\nc\bepsilon{{\mbox{\boldmath $\epsilon$}}}
  \nc\bomegas{{\mbox{\boldmath\scriptsize $\omega$}}}
  \nc\bepsilons{{\mbox{\boldmath \scriptsize$\epsilon$}}}

  \nc\bxi{{\mbox{\boldmath $\xi$}}}
\nc\bmu{{\mbox{\boldmath $\mu$}}} \nc\bcN{{\mbox{\boldmath $\cal{N}$}}} \nc\bcm{{\mbox{\boldmath $\cal{M}$}}} \nc\blambda{{\mbox{\boldmath
$\lambda$}}}

\newcommand{\Tmn}{\bold{T}_{\lambda^1, \lambda^2}^{\nu}}

\newcommand{\lie}[1]{\mathfrak{#1}}
\newcommand{\ol}[1]{\overline{#1}}
\makeatletter
\def\section{\def\@secnumfont{\mdseries}\@startsection{section}{1}%
  \z@{.7\linespacing\@plus\linespacing}{.5\linespacing}%
  {\normalfont\scshape\centering}}
\def\subsection{\def\@secnumfont{\bfseries}\@startsection{subsection}{2}%
  {\parindent}{.5\linespacing\@plus.7\linespacing}{-.5em}%
  {\normalfont\bfseries}}
\makeatother
\def\subl#1{\subsection{}\label{#1}}
 \nc{\Hom}{\operatorname{Hom}}
  \nc{\mode}{\operatorname{mod}}
\nc{\End}{\operatorname{End}} \nc{\wh}[1]{\widehat{#1}} \nc{\Ext}{\operatorname{Ext}}
 \nc{\ch}{\operatorname{ch}} \nc{\ev}{\operatorname{ev}}
\nc{\Ob}{\operatorname{Ob}} \nc{\soc}{\operatorname{soc}} \nc{\rad}{\operatorname{rad}} \nc{\head}{\operatorname{head}}
\def\Im{\operatorname{Im}}
\def\gr{\operatorname{gr}}
\def\mult{\operatorname{mult}}
\def\Max{\operatorname{Max}}
\def\ann{\operatorname{Ann}}
\def\sym{\operatorname{sym}}
\def\loc{\operatorname{loc}}
\def\Res{\operatorname{\br^\lambda_A}}
\def\und{\underline}
\def\Lietg{$A_k(\lie{g})(\bs_\xiigma,r)$}
\def\res{\operatorname{res}}

 \nc{\Cal}{\cal} \nc{\Xp}[1]{X^+(#1)} \nc{\Xm}[1]{X^-(#1)}
\nc{\on}{\operatorname} \nc{\Z}{{\bold Z}} \nc{\J}{{\cal J}} \nc{\C}{{\bold C}} \nc{\Q}{{\bold Q}}
\renewcommand{\P}{{\cal P}}
\nc{\N}{{\Bbb N}} \nc\boa{\bold a} \nc\bob{\bold b} \nc\boc{\bold c} \nc\bod{\bold d} \nc\boe{\bold e} \nc\bof{\bold f} \nc\bog{\bold g}
\nc\boh{\bold h} \nc\boi{\bold i} \nc\boj{\bold j} \nc\bok{\bold k} \nc\bol{\bold l} \nc\bom{\bold m} \nc\bon{\bold n} \nc\boo{\bold o}
\nc\bop{\bold p} \nc\boq{\bold q} \nc\bor{\bold r} \nc\bos{\bold s} \nc\boT{\bold t} \nc\boF{\bold F} \nc\bou{\bold u} \nc\bov{\bold v}
\nc\bow{\bold w} \nc\boz{\bold z} \nc\boy{\bold y} \nc\ba{\bold A} \nc\bb{\bold B} \nc\bc{\mathbb C} \nc\bd{\bold D} \nc\be{\bold E} \nc\bg{\bold
G} \nc\bh{\bold H} \nc\bi{\bold I} \nc\bj{\bold J} \nc\bk{\bold K} \nc\bl{\bold L} \nc\bm{\bold M}  \nc\bo{\bold O} \nc\bp{\bold
P} \nc\bq{\bold Q} \nc\br{\bold R} \nc\bs{\bold S} \nc\bt{\bold T} \nc\bu{\bold U} \nc\bv{\bold V} \nc\bw{\bold W} \nc\bx{\bold
x} \nc\KR{\bold{KR}} \nc\rk{\bold{rk}} \nc\het{\text{ht }}
\nc\bz{\mathbb Z}
\nc\bn{\mathbb N}

\nc\toa{\tilde a} \nc\tob{\tilde b} \nc\toc{\tilde c} \nc\tod{\tilde d} \nc\toe{\tilde e} \nc\tof{\tilde f} \nc\tog{\tilde g} \nc\toh{\tilde h}
\nc\toi{\tilde i} \nc\toj{\tilde j} \nc\tok{\tilde k} \nc\tol{\tilde l} \nc\tom{\tilde m} \nc\ton{\tilde n} \nc\too{\tilde o} \nc\toq{\tilde q}
\nc\tor{\tilde r} \nc\tos{\tilde s} \nc\toT{\tilde t} \nc\tou{\tilde u} \nc\tov{\tilde v} \nc\tow{\tilde w} \nc\toz{\tilde z} \nc\woi{w_{\omega_i}}
\nc\chara{\operatorname{Char}}
\title[Macdonald Polynomials and level two Demazure modules ]
{ Macdonald Polynomials and  level two Demazure modules for affine $\lie{sl}_{n+1}$ }
\author[Biswal, Chari, Shereen, Wand]{Rekha Biswal, Vyjayanthi Chari, Peri Shereen, Jeffrey Wand}
\address{Max-Planck-Institut fur Mathematik\\
  Vivatsgasse 7, 53111 Bonn, Germany } 
 \email{rekha@mpim-bonn.mpg.de}
 \address{Department of Mathematics\\ 
  University of California, Riverside\\ 
  900 University Ave, Riverside, CA 92521} \thanks{VC was partially supported by DMS- 1719357, the Simons Fellows Program and by the  Infosys Visiting Chair position at the Indian Institute of Science.}
\email{chari@math.ucr.edu}
\address{Department of Mathematics and Statistics\\ California State University, Monterey Bay\\ 100 Campus Center, Seaside, CA 93955}\email{pshereen@csumb.edu}
\address{Department of Mathematics and Statistics\\ California State University, Monterey Bay\\ 100 Campus Center, Seaside, CA 93955}\email{jwand@csumb.edu}
\begin{abstract}
 We define  a family of  symmetric polynomials $G_{\nu,\lambda}(z_1,\cdots, z_{n+1},q)$  indexed by a pair of  dominant integral weights for a root system of type $A_n$. The polynomial  $G_{\nu,0}(z,q)$ is  the specialized Macdonald polynomial $P_\nu(z,q,0)$ and is known to be the graded character of a level one Demazure module associated to the affine Lie algebra $\wh{\lie{sl}}_{n+1}$.  We prove that    $G_{0,\lambda}(z,q)$ is the   graded character of a level two Demazure module for $\wh{\lie{sl}}_{n+1}$. Under suitable conditions on $(\nu,\lambda)$ (which apply to the pairs $(\nu,0)$ and $(0,\lambda)$) we  prove that $G_{\nu,\lambda}(z,q)$ is Schur positive, i.e., it 
  can be written as a linear combination of  Schur polynomials with coefficients in $\mathbb Z_+[q]$.
 We further prove that  $P_\nu(z,q,0)$ 
is  a  linear combination of  elements $G_{0,\lambda}(z,q)$ with  the coefficients being essentially products of $q$-binomials.
Together with a result of K. Naoi, a consequence of our result is  an explicit  formula for the  specialized Macdonald polynomial associated to   a non-simply laced  Lie algebra  as a linear combination of level one Demazure characters.
 
 \end{abstract}
 \maketitle
 
 \section*{Introduction}
 In 1987, I.G.Macdonald introduced a family of orthogonal symmetric polynomials $P_\lambda(z,q,t)$, $z=(z_1,\cdots, z_n)$  which are  a basis for the ring of of symmetric polynomials in $\bc(q,t)[z_1,\cdots, z_n]$; here $\lambda$ varies over the set of partitions of length at most $n$. These polynomials interpolate between several well--known families of symmetric polynomials such as the  Schur polynomials $P_\lambda(z,0,0)$, the Hall-Littlewood polynomials $P_\lambda(z,0,t)$ and the Jack polynomials to name a few. The subject has deep connections with  combinatorics, geometry and representation theory and there is a vast literature on the subject.
 
 In this paper we shall be interested in the connection between the specialized Macdonald polynomials $P_\lambda(z,q,0)$ and the representation theory of the  affine Lie algebra $\wh{\lie{sl}}_{n+1}$. Such a connection was first shown to exist in \cite{San} and we  discuss this briefly. Fix a Borel subalgebra $\widehat{\lie b}$ of $\wh{\lie{sl}}_{n+1}$, let $\Lambda_i$ $0\le i\le n$ be a corresponding set of fundamental weights and let $V(\Lambda_i)$ be the associated  integrable highest weight representation of $\wh{\lie{sl}}_{n+1}$. Given a partition $\lambda$ or equivalently a dominant integral weight of $\lie{sl}_{n+1}$ there exists an element $w$  of the affine Weyl group and $0\le i\le n$ such that $w_\circ
w\Lambda_i=(\lambda+\Lambda_0)$ where $w_\circ$ is the longest element of the Weyl group of $\lie{sl}_{n+1}$. The main result of \cite{San} shows that  $P_\lambda(z,q,0)$ is the character of  the $\widehat{\lie b}$--module  generated by the one-dimensional subspace $V_{w\Lambda_i}(\Lambda_i)$ of $V(\Lambda_i)$. This family of Demazure modules is special since the   modules admit an action of  the standard maximal parabolic subalgebra.  (We remind the reader that the derived subalgebra of the standard maximal parabolic is also called  the current algebra of $\lie{sl}_{n+1}$; i.e., it is the subalgebra of $\wh{\lie{sl}}_{n+1}$   consisting of polynomial maps $\bc\to \lie{sl}_{n+1}$.)  This result was later extended in  \cite{Ion} 
to twisted affine Lie algebras and the  untwisted affine Lie algebras associated to a simply-laced simple Lie algebra. The result was  known to be false for the affine algebras associated to the non-simply laced Lie algebras.

Recently it was shown in \cite{CI} that  $P_\lambda(z,q,0)$ can in fact be realized as the character of a suitable module for the current algebra associated to a non-simply laced simple Lie algebra. This module has  the corresponding  Demazure module as a (possibly) proper quotient. However, the character of the Demazure module is not known in these cases and   one of the goals of this paper is to solve this problem. One  approach to the problem is to use the idea of  Demazure flags which was first introduced and developed in \cite{J}. It was further studied in  \cite{Naoi}  where it is  shown that in the case of $B_r, C_r, F_4$ it suffices to study the relationship between level one and level two Demazure modules for affine $\lie{sl}_{n+1}$. (In the case of $G_2$ one also has to understand level three modules for $A_1$; this was done in \cite{BCSV}, \cite{CSSW},  and we will say no more about it in this paper). The level two modules  are also modules for the current algebra, are indexed by dominant integral weights $\lambda$ and are $\widehat{\lie b}$-submodules   of $V(\Lambda_i+\Lambda_j)$; here  $0\le i,j\le n$ and  $w$ are chosen so that $w_\circ w(\Lambda_i+\Lambda_j)=(\lambda+2\Lambda_0)$.  As a consequence of the main result of this paper, we give an explicit formula for the specialized  Macdonald polynomials associated to root systems of type $B,C,F$ in terms of the level one affine Demazure modules for the associated affine Lie algebra. 

The level two Demazure modules for affine $\lie{sl}_{n+1}$ also appear in a completely different context. They are the classical limit of an important family of modules for quantum affine $\lie{sl}_{n+1}$. These modules occur in the work of \cite{HL1}, \cite{HL2}  and we  refer the reader to Section \ref{motivation} of this paper for further details. Not much is known about the structure of these modules for the quantum affine algebra. Our results show that their character is given by the polynomials $G_{0,\lambda}(z,q)$.

To understand the characters of level two Demazure modules it is best to work in a more general framework.  Thus we define a family of (finite-dimensional) modules $M(\nu,\lambda)$ for the current algebra of $\lie{sl}_{n+1}$  which are indexed by two dominant integral weights and interpolate between level one and level two Demazure modules. The graded characters of these modules are precisely the polynomials $G_{\nu,\lambda}(z,q)$. Although in this paper we only consider the case of $\lie{sl}_{n+1}$ and level two, the definitions we give, go through in a straightforward way to other simple Lie algebras and higher levels. But, the representation theory becomes much more difficult. The crucial input in this paper coming from the connections with quantum affine algebra associated to $\lie{sl}_{n+1}$ is 
 missing in the other cases for the moment (see however \cite{CDM}).
 
 The paper is organized as follows. In Section \ref{sec1} we define the polynomials $G_{\nu,\lambda}(z,q)$ (see equations  \eqref{defg} and \eqref{gnl})  in terms of Macdonald polynomials (of type $A_n$) and state the main theorem. We  explain in some detail in Section \ref{motivation} the motivation and the methods used to prove the main results  and discuss possible further directions. In Section \ref{mnl} we introduce the modules $M(\nu,\lambda)$ and study their graded characters. In Section \ref{strategy} we explain the key steps in deducing the main theorem of this paper.
 In the remaining two sections of the paper
we prove the key step and   give explicit formulae for the specialized Macdonald polynomials in term of the polynomials $G_{0,\lambda}$ and vice-versa.

\medskip

{\em{Acknowledgements. The second author thanks Matheus Brito for many helpful discussions. She is grateful to Bogdan Ion for comments on a preliminary version. The authors thank S. Viswanath for several comments which improved the overall exposition of the paper.}}
\section{The main results}\label{sec1}
 Let $\bc$, $\bz$, $\bz_+$ and  $\bn$ be the set of complex numbers,  integers, non--negative integers and positive  integers respectively. For $i,j\in \bz$ with $i\le j$ we let $[i,j]=\{i,i+1,\cdots, j\}$. Given an indeterminate $q$ and $n,r\in\bz$ set \begin{gather*}[n]_q=\frac{1-q^n}{1-q},\ \ \qbinom{n}{r}_q=\frac{[n]_q[n-1]_q\cdots[n-r+1]_q}{[1]_q[2]_q\cdots[r]_q},\ \ \  0<r\le n,\\ \qbinom{n}{0}_q=1,\ \ n\in\bz_+, \ \ \qbinom{n}{r}_q=0 \ \  {\rm{unless}} \ \ \{n,r,n -r\}\subset\bz_+.\end{gather*}   

\subsection{}
Let $\lie g$ be the simple Lie algebra $\lie{sl}_{n+1}$  and let $\lie h$ be a fixed Cartan subalgebra. Fix a set  $\{\alpha_i: 1\le i\le n\}$ of simple roots for $(\lie g,\lie h)$  and  a corresponding  set $\{\omega_i: 1\le i\le n\}$  of fundamental weights. It is  convenient to set $\omega_0=\omega_{n+1}=0$. Let $$R^+=\{\alpha_{i,j}=\alpha_i+\alpha_{i+1}+\cdots+\alpha_j:1\le i\le j\le n\},$$  be the set of positive roots. Let $Q$, $P$ (resp. $Q^+,\  P^+$) be  the $\bz$-span (resp. $\bz_+$-span)  of the set of  simple roots and fundamental weights respectively.    
Define a partial order $\le $  on $P$  by: $\mu\le\lambda$ iff $\lambda-\mu\in Q^+$.
Given $\lambda=\sum_{i=1}^n r_i\omega_i\in P^+$ and $\eta=\sum_{i=1}^nr_i\alpha_i\in Q^+$  set $$\Ht\lambda=\sum_{i=1}^nr_i=\Ht_r\eta.$$
Let $(\ ,\ ): P\times P\to \mathbb Q $ be the symmetric bilinear form  defined by setting $(\omega_i,\alpha_j)=\delta_{i,j}$ and set  $$ P^+(1)=\{\lambda\in P^+: (\lambda,\alpha_i)\le 1\ \  {\rm{for\ all}}\  \ 1\le i\le n\}.$$ 
{\em In  the rest of the paper we shall write   (often without further mention) an element $\lambda\in P^+$ as  a sum $\lambda=2\lambda_0+\lambda_1$ with $\lambda_0\in P^+ $ and $\lambda_1\in P^+(1)$. } 
\subsection{The polynomials $p_\lambda^\mu(q)$ and  $G_{\nu,\lambda}(z,q)$} 
Given an indeterminate $q$
and  $
\lambda,\mu\in P^+$   define  $p_\lambda^\mu\in\mathbb Z_+[q]$ by 
 $$p_\lambda^\mu(q) =q^{\frac12(\lambda+\mu_1,\  \lambda-\mu)}\displaystyle\prod_{j=1}^n\qbinom{(\lambda-\mu,\ \omega_j) + (\mu_0, \alpha_j)}{(\lambda-\mu,\omega_j)}_q,\ \ \mu=2\mu_0+\mu_1.$$
 Notice that $p_\lambda^\lambda=1$ and  $p_\lambda^\mu=0$ if $\lambda-\mu\notin Q^+$. Moreover, $$ (\lambda+\mu_1,\lambda-\mu)
=(\lambda-\mu,\lambda-\mu) +2(\mu-\mu_0,\lambda-\mu)\in 2\bz_+,\ \ {\rm{if}}\ \ \lambda-\mu\in Q^+,$$ and in particular  $p_\lambda^\mu\in\bz_+[q]$ as asserted. It will be convenient to extend the definition by setting $p_\lambda^\mu=0$ if   $\lambda$ or $\mu$ are in $P\setminus P^+$.
\medskip

\noindent Given  $\lambda\in P^+$ and a set of indeterminates $z=(z_1,\cdots, z_{n+1})$ let $s_\lambda(z)$ and $P_\lambda (z,q,0)$ be the  associated Schur polynomial and the  specialized Macdonald polynomial  respectively.  Recall that the Schur  polynomials  are a basis for the ring of symmetric polynomials  $\mathbb C[z_1,\cdots, z_{n+1}]$ while the specialized Macdonald polynomials are a basis for the ring of symmetric polynomials in  $\mathbb C(q)[z_1,\cdots, z_{n+1}]$.

\medskip

\noindent Define  elements  $G_\lambda(z,q)\in\bc[q][z_1,\cdots,z_{n+1}]$ recursively
  by requiring,
  \begin{gather}
  \label{defg} P_\lambda(z,q,0)= \sum_{{\mu\in P^+}} p_\lambda^\mu(q)G_\mu(z,q).\end{gather}
 Since $p_\lambda^\mu=0$ if $\mu\nleq \lambda$, it follows that \begin{gather*}P_{\omega_i}(z,q,0)= G_{\omega_i}(z,q)=s_{\omega_i}(z),\ \ i\in [0,n],\ \
 P_\lambda(z,q,0)= G_\lambda(z,q)+\sum_{\mu<\lambda}p_\lambda^\mu G_\mu(z,q).\end{gather*}
  Hence  the elements $\{G_\lambda(z,q):\lambda\in P^+\}$ are also   a basis for the ring of symmetric polynomials in   $\mathbb C(q)[z_1,\cdots, z_{n+1}]$.  Hence  there exists  a subset  $\{a_\lambda^\mu: \mu\in P^+\}\subset \bc(q)$  such that $$a_\lambda^\lambda=1,\ \ a_\lambda^\mu=0\ \ {\rm{if}}\ \  \lambda-\mu\notin Q^+,$$ satisfying \begin{equation}\label{defg}G_\lambda(z,q)=\sum_{\mu\in P^+}a_\lambda^\mu(q)P_{\mu}(z,q,0)  \ \ \ {\rm{and}} \ \ 
 \ \ \ \sum_{\mu\in P^+}a_\lambda^\mu p_\mu^\nu=\delta_{\lambda,\mu}=\sum_{\mu\in P^+}p_\lambda^\mu a_\mu^\nu.\end{equation}

\noindent  Given a pair of elements $\nu,\lambda\in P^+$ set
\begin{equation}\label{gnl}G_{\nu,\lambda}(z,q)=\sum_{{\mu\in P^+}}q^{(\lambda+\nu-\mu,\nu)}a_\lambda^{\mu-\nu}(q)P_{\mu}(z,q,0),\end{equation} where we understand that $a_\lambda^{\mu-\nu}=0$ if $\mu-\nu\notin P^+$.
 Notice that $$G_{\nu,0}(z,q)= P_\nu(z,q,0),\ \ \  G_{0,\lambda}(z,q)= G_\lambda(z,q).$$ Here the second equality is obvious while the  first equality follows since $a_0^0=1$ and $a_0^{\mu-\nu}\ne 0$ only if $\mu-\nu\in P^+$ and  $\nu-\mu\in Q^+$ which is possible only if $\nu=\mu$. 
\subsection{The main result} \label{mainresult} For $\mu\in P^+$ set $$\max\mu=\max\{i:\mu(h_i)>0\},\ \ \min\mu=\min\{i:\mu(h_i)>0\}.$$ We say that  a pair $(\nu,\lambda)\in P^+\times P^+$ is admissible if one of the following hold: writing $\lambda=2\lambda_0+\lambda_1$, $\nu=2\nu_0+\nu_1$, either,\\ $\bullet$ $\lambda_1=0$, or\\ $\bullet$ $\lambda_1\ne 0$,  $\nu_0=\omega_i$ for some $i\in[0,n]$,  $ \max\nu_1<\min\lambda_1$  and if $  i\in[1,n]$ then  $ i<\min\lambda_1-1$ and $   \nu_1(h_i)=\nu_1(h_{i+1})=0.$
\medskip

A symmetric polynomial in $\bc(q)[z_1,\cdots, z_{n+1}]$ is Schur-positive if it is in  the $\bz_+[q]$ span of  $\{s_\mu: \mu\in P^+\}$.
  The following is one of the main results of this paper. 
 \begin{thm}\label{main1a} For an admissible pair $(\nu,\lambda)\in P^+\times P^+$  the polynomial $G_{\nu,\lambda}(z,q)$ is Schur positive. Moreover for such pairs we have $$G_{\nu,\lambda}(z,1)= G_{\nu,0}(z,1)G_{0,\lambda}(z,1).$$
 \end{thm}
 
 \medskip
 
\begin{rem} We also give  a closed formula for the polynomials $a_\lambda^\mu$.   Unlike the polynomial $p_\lambda^\mu$ which is non--zero for all $\mu\in P^+$ with $\mu\le \lambda$, the polynomials $a_\lambda^\mu$ are non--zero on a much smaller set and the statement requires some additional definitions;  the formulae can be found in    Section \ref{secfour} (see also Section \ref{strategy} and equations \eqref{qbinom1}, \eqref{geng}). \end{rem}
\subsection{Examples}  
 We give some explicit examples of  $G_{\nu,\lambda}(z,q)$ in terms of Schur polynomials.
 
Suppose that  $\lie g=\lie{sl}_3$. Then,
\begin{eqnarray*}
&G_{0,\omega_1+\omega_2}&=s_{\omega_1+\omega_2},\\
&G_{0,2\omega_1+2\omega_2}&=s_{2\omega_1+2\omega_2}+qs_{\omega_1+\omega_2}+q^2 s_{0},\\
&G_{0,2\omega_1+3\omega_2}&=s_{2\omega_1+3\omega_2}+q[2]_qs_{\omega_1+2\omega_2}+q^2s_{2\omega_1}+q^2[2]_qs_{\omega_2}+qs_{3\omega_1+\omega_2},\\
&G_{0,3\omega_1+3\omega_2}&=s_{3\omega_1+3\omega_2}+q(s_{4\omega_1+\omega_2}+s_{\omega_1+4\omega_2})+(2q^2+q)s_{2\omega_1+2\omega_2}\\&&+q^2[2]_q(s_{3\omega_1}+s_{3\omega_2})+(q[2]_q)^2s_{\omega_1+\omega_2}+q^3s_0,\\
&G_{\omega_1+\omega_2,2\omega_2}&=s_{\omega_1+3\omega_2}+qs_{2\omega_1+\omega_2}+ q^2s_{\omega_1}+q[2]_qs_{2\omega_2}.
\end{eqnarray*}

Suppose that  $\lie g=\lie{sl}_4$. Then,
\begin{eqnarray*}&G_{0,\omega_1+\omega_2+\omega_3}&=s_{\omega_1+\omega_2+\omega_3}+qs_{\omega_2},\\
&G_{0,2\omega_1+\omega_2+\omega_3}&=s_{2\omega_1+\omega_2+\omega_3}+q(s_{\omega_1+2\omega_3}+s_{\omega_1+\omega_2})+q^2s_{\omega_3}.
\end{eqnarray*}

Suppose that  $\lie g=\lie{sl}_5$. Then,
\begin{eqnarray*}&G_{0,\omega_1+\omega_4}&=s_{\omega_1+\omega_4},\\
&G_{0,2\omega_1+2\omega_4}&=s_{2\omega_1+2\omega_4} +qs_{\omega_1+\omega_4} + q^2s_0,\\
&G_{0,\omega_1+\omega_2+\omega_4}&=s_{\omega_1+\omega_2+\omega_4}+qs_{\omega_2}.
\end{eqnarray*}

It seems likely that the polynomials $G_{\nu,\lambda}$ are always Schur positive. We give some examples where this holds  for non-admissible pairs $(\nu,\lambda)$.
In the case of  $\lie{sl_3}$ we have $$G_{2\omega_1,\omega_1+\omega_2}=s_{3\omega_1+\omega_2}+q[2]_q(s_{\omega_1+2\omega_2}+s_{2\omega_1}+qs_{\omega_2}).$$
In the case of $\lie{sl_4}$ we have  
\begin{eqnarray*}&G_{\omega_3,\omega_1+\omega_2+\omega_3}&=s_{\omega_1+\omega_2+2\omega_3}+q(s_{2\omega_1+\omega_3}+s_{\omega_1+2\omega_2})+q^2 s_{\omega_1}+q[2]_qs_{\omega_2+\omega_3},\\
  &G_{\omega_2,\omega_1+\omega_2+\omega_3}&=s_{\omega_1+2\omega_2+\omega_3}+q(s_{2\omega_1+2\omega_3}+s_{2\omega_1+\omega_2}+s_{\omega_2+2\omega_3})+q[2]_qs_{2\omega_2}+2q^2s_{\omega_1+\omega_3}+q^3s_0,\\
&G_{\omega_1,\omega_1+\omega_2+\omega_3}&=s_{2\omega_1+\omega_2+\omega_3}+q(s_{2\omega_2+\omega_3}+s_{\omega_1+2\omega_3})+q[2]_qs_{\omega_1+\omega_2}+q^2s_{\omega_3}.\end{eqnarray*}
For  $\lie{sl}_5$, we have
\begin{eqnarray*}
&G_{\omega_1,\omega_1+\omega_4}&=s_{2\omega_1+\omega_4}+q(s_{\omega_2+\omega_4} + s_{\omega_1}),\\
&G_{\omega_1+\omega_2,\omega_1+\omega_4}&=s_{2\omega_1+\omega_2+\omega_4}+q(s_{3\omega_1}+s_{2\omega_1+\omega_4})+q[2]_q(s_{\omega_1+\omega_3+\omega_4}+qs_{\omega_3})\\
&&+q^3s_{2\omega_4}+(2q^2+q)s_{\omega_1+\omega_2}.
\end{eqnarray*}

 \subsection{}\label{motivation}  We explain in some detail the motivation for these results and  assume for just this discussion, that  $\lie g$ is an arbitrary simple Lie algebra. As discussed  in the introduction, it was shown in  \cite{Ion, San} that  the specialized  Macdonald polynomial $P_\lambda(z,q,0)$ is the character of a Demazure module occurring in a fundamental integrable highest weight representation of the affine Lie algebra associated to a simple Lie algebra of type $A,D,E$.  The graded characters of Demazure modules occurring in higher level integrable representations are not very well  understood. There are however some combinatorial  results for weights  of the form $r\lambda$ for some $r\in\bz_+$ and $\lambda$ a dominant integral weight for $\lie g$ (see for instance \cite{L,LNSSS, LS}).  However, in the case of $\lie{sl}_{n+1}$ for instance, it is the study of  level two  Demazure modules  which are {\em not} of this  type that are  of  interest. This interest  arises from the  connections with cluster algebras and monodial categorification (see  \cite{BC, HL1, HL2}). It was shown in \cite{BCM} that the classical limit $q\to 1$ of the representations corresponding to a cluster variable in a cluster algebra of type $A_n$ is a level two Demazure module.
  This is one motivation for our study since our results show that the  character of the irreducible representation of the quantum affine algebra is given by the polynomials $G_{0,\lambda}(q,z)$, $\lambda\in P^+(1)$.
 
 Another reason to be  interested in this problem comes from  the connection with Macdonald polynomials associated to root systems of simple Lie algebras of type $B,C,F,G$. For these algebras it has long been known that the Macdonald polynomial is \lq\lq too big\rq\rq \ to be the character of the Demazure module. It was shown recently in \cite{CI} that for these algebras the  Macdonald polynomial is the character of a family of modules called the local Weyl modules. These modules denoted $W_{\loc}(\lambda)$, are defined for all simple Lie algebras and  are indexed by the dominant integral weights of the underlying simple Lie algebra. Like the Demazure modules they are modules for the standard maximal parabolic subalgebra of the affine Lie algebra or the current algebra. In types $A,D, E$ the Weyl and Demazure modules coincide (see \cite{CL}, \cite{FL}) but in other types the corresponding Demazure module can be  a proper quotient of the local Weyl module. 
 
 A  study of the relationship between the local Weyl modules and the level one Demazure modules in the non-simply laced case  was initiated  by K. Naoi in \cite{Naoi} and we now discuss his results.
Let $\mu$ be a dominant integral weight for a simple Lie algebra of type $X_n$ and  $D(1,\mu)$    a
Demazure module occurring in a level one integrable highest weight representation of the affine Lie algebra of type $X_n^{(1)}$ where $X\in\{B,C, F\}$.  Naoi proved that if $\mu$ is zero on the short simple roots then the local Weyl module is isomorphic to a level one Demazure module. 
Otherwise, he showed that  the local Weyl module admits a non-trivial flag whose successive quotients are isomorphic to level one Demazure modules. The modules in question are all graded and hence one can define the graded multiplicity of the  module $D(1,\mu)$ in  $W_{\loc}(\lambda)$, denoted $[W_{\loc}(\lambda):D(1,\mu)]$. 
 
 Naoi's next result related this graded multiplicity to a problem  in type $A$.  Let $\mu_s\in P$ be the element defined by $\mu_s(h_i)=0$ if $\alpha_i$ is a long roots and $\mu_s(h_i)=\mu(h_i)$ if $\alpha_i$ is short. We also regard $\mu_s$ as  a dominant integral weight for the simple Lie algebra $\lie g_s$ which is   generated by the short simple roots. Note that $\lie g_s$ is of type $A_1$ if $X_n=B_n$, of type $A_2$ if  $X_4=F_4$ and of type $A_{n-1}$ if $X_n=C_n$.
Let $W^s_{\loc}(\mu_s)$ and $D^s(2,\mu_s)$ be the local Weyl module  (equivalently the level one Demazure module) and the level two Demazure module associated to the  untwisted affine Lie algebra $\wh{\lie g}_s$.  Naoi proved that the module $W^s_{\loc}(\lambda_s)$ admits a flag whose successive quotients are level two Demazure modules $D^s(2,\mu_s)$ and moreover,
\begin{gather*}[W_{\loc}(\lambda): D(1,\mu)]=\delta_{\lambda-\lambda_s, \mu-\mu_s}[W^s_{\loc}(\lambda_s): D^s(2,\mu_s)].\end{gather*} 
Together with Naoi's results, Theorem \ref{main1a} of this paper can be reformulated as asserting the following: if $\lie g$ is of type $B,C,F$ then 

\begin{gather*}
P_\lambda(z,q,0)=\sum_{\mu\in P^+} \delta_{\mu-\mu_s,\lambda-\lambda_s}p_{\lambda_s}^{\mu_s} \ch_{\gr} D(1,\mu),\ .\end{gather*}
In particular our result gives a   formula for the graded character of the level one Demazure module when $\lie g$ is not simply laced in terms of Macdonald polynomials. This finally completes the picture begun in  \cite{Ion, Ion2, San}.

\subsection{} The proof of our results is somewhat indirect.
We introduce a family of graded finite--dimensional  modules $M(\nu,\lambda)$ for the  standard maximal parabolic subalgebra of affine $\lie{sl}_{n+1}$. These  interpolate between the local Weyl module $M(\nu,0)=W_{\loc}(\nu)$  and the level  two Demazure module $M(0,\lambda)= D(2,\lambda)$. The graded characters of the modules $M(\nu,0)$, $\nu\in P^+$ (resp. $M(0,\lambda)$, $\lambda\in P^+$) are linearly independent and their $\bz[q]$-span contains the graded character of $M(\nu,\lambda)$. 
We show (see Proposition \ref{ses} and Lemma \ref{sescon}) that if $(\nu,\lambda)$ is admissible then the graded characters of $M(\nu,\lambda)$ satisfy certain recursions.  These recursions can be solved  to express the character of $M(\nu,\lambda)$ in terms of $M(\mu,0)$ (resp. $M(0,\mu)$). Together with the results in \cite{San} this leads to a proof of Theorem \ref{main1a}.  

\subsection{} There are some obvious questions that can be raised. For instance, one could define the polynomials $G_{\nu,\lambda}$ in the same way for other simple Lie algebras and ask if they arise as the graded character of some module. Another question would be to ask if there is some general framework in which to study polynomials which interpolate between characters of  higher level Demazure modules. In fact, one can define suitable analogs of the modules $M(\nu,\mu)$ to address these questions. But the difficulty in answering these questions lies in the further  development of  the appropriate representation theory. In the particular case addressed in this paper,  the input coming from the representation theory of quantum affine $\lie{sl}_{n+1}$ \cite{BC, BCM, HL1, HL2} has been critical. But the analogs of these results are not known and are not easily formulated for the  other quantum affine algebras or for higher level Demazure modules (see however \cite{CDM}). We hope to return to these ideas elsewhere.

\section{The modules $M(\nu,\lambda)$} \label{mnl} We recall the definition of the current algebra associated to $\lie{sl}_{n+1}$. We also  recall 
 several results on their representation theory which are needed for our study. We  then  introduce the modules $M(\nu,\lambda)$ and  explain in Section \ref{strategy} 
 how the study of these modules leads to   the proof of Theorem \ref{main1a}. In the rest of the section we study the graded characters of these modules and show that they satisfy certain recursions.

\subsection{The current algebra $\lie{sl}_{n+1}[t]$}
Let $t$ be an indeterminate and   $\bc[t]$ the  polynomial ring with complex coefficients. 
Denote by  $ \lie g[t]$ the Lie algebra  with underlying vector space $ \lie g\otimes \bc[t]$ and commutator given by $$[a\otimes f, b\otimes g]=[a,b]\otimes fg,\ \ a,b\in \lie g,\ \ f,g\in\bc[t].$$ Then $ \lie g[t]$ and its universal enveloping algebra admit a natural $\bz_+$-grading given by declaring a monomial $(a_1\otimes t^{r_1})\cdots (a_p\otimes t^{r_p})$ to have grade $r_1+\cdots +r_p$, where $a_s\in \lie g$ and $r_s\in\bz_+$ for $1\le s\le p$. We  freely identify the subspace $\lie g\otimes 1$ with the Lie algebra $\lie g$.
 
 We shall be interested in the category of finite--dimensional $\bz$--graded modules for $\lie g[t]$. An object of this category is a  finite-dimensional module $V$ for $ \lie g[t]$ with a compatible  $\bz$--grading,  $$V=\bigoplus_{s\in\bz} V[s],\ \ (x\otimes t^r)V[s]\subset V[r+s],\ \ x\in \lie g,\ \ r\in\bz_+.$$ For $s\in\bz$, clearly $V[s]$ is a $\lie g$--module.  For any $p\in\bz$  let $\tau_p^*V$ be the graded $\lie g[t]$--module which is given by shifting all the grades up by $p$ and leaving the action of $ \lie g[t]$ unchanged. The morphisms between graded modules are $ \lie g[t]$--module maps  of grade zero.

 Let $\bz[q,q^{-1}][P]$ be the group ring of $P$ with basis $\{e(\mu):\mu\in P\}$ and coefficients in $\bz[q,q^{-1}]$. Any finite-dimensional $\lie g$-module  $V$  can be written as, 
 
 $$V=\bigoplus_{\mu\in P} V_\mu,\ \  V_\mu=\{v\in V: hv=\mu(h)v,\ \ h\in\lie h\},$$
 and the character of a $\lie g$--module $V$ is the element $\ch V=\sum_{\mu\in P}\dim V_\mu\  e(\mu)$ of $\bz[P]$. 
 The irreducible finite--dimensional $\lie g$-modules are indexed by elements of $P^+$, and  given $\lambda\in P^+$ we let $V(\lambda)$ be an irreducible module corresponding to $\lambda$.
 Setting  $z_i = e(\omega_i-\omega_{i-1})$, $1\le i\le n+1$, a classical result  asserts that $\ch V(\lambda)$ is just the Schur polynomial $s_\lambda(z)$.
  If  $V$ is a graded $\lie g[t]$--module, let $$\ch_{\gr} V=\sum_{s\in\bz} q^s\ch V[s]=\sum_{\mu\in P^+}\left( \sum_{s\in\bz} \dim\Hom_{\lie g}(V(\mu), V[s])q^s\right)\ch V(\mu) .$$ 
 
 \subsection{The modules $M(\nu,\lambda)$}  We recall the definition of a  family of modules $M(\nu,\lambda)$ which were  introduced and studied in \cite{Wand}.  Fix a Chevalley basis $\{x_{i,j}^\pm, h_i: 1\le i\le j\le n\}$  for $\lie g$ and let $x_{j,j}^\pm=x_j^\pm$, $h_{i,j}=h_i+\cdots+h_j$ for all $1\le i\le j\le n$.

For $\nu,\lambda\in P^+$ with  $\lambda=2\lambda_0+\lambda_1$, let $M(\nu,\lambda)$ be the $\lie g[t]$--module generated by an element $w_{\nu,\lambda}$ with the following relations: 
\begin{equation}\label{localweyld} (x_i^+\otimes 1)w_{\nu,\lambda}=0,\ \ (h\otimes t^r)w_{\nu,\lambda}=\delta_{r,0}(\lambda+\nu)(h)w_{\nu,\lambda},\ \ \  (x_i^-\otimes 1)^{(\lambda+\nu)(h_i)+1}w_{\nu,\lambda}=0,\end{equation} 
 \begin{equation}\label{lnu}(x_\alpha^-\otimes t^{\nu(h_\alpha)+\lceil\lambda(h_\alpha)/2\rceil})w_{\nu,\lambda}=0,\end{equation} for all $i\in[1,n]$, $h\in\lie h$, $r\in\bz_+$  and $\alpha\in R^+$.
An inspection of the defining relations shows that 
\begin{gather}\label{fund} M(\omega_i,0)\cong_{\lie g[t]} M(0,\omega_i)\cong_{\lie g} V(\omega_i),\ \ M(0,2\omega_i)\cong_{\lie g} V(2\omega_i),\ \  i\in[1,n]\\
\label{iso1} M(\nu+\omega_i, 2\lambda_0)\cong_{\lie g[t]} M(\nu, 2\lambda_0+\omega_i),\ \\\label{iso2} M(\nu,\omega_i)\cong_{\lie g[t]} M(\nu+\omega_i,0).\end{gather} 
Since the defining relations of $M(\nu,\lambda)$ are graded, it follows that the module $M(\nu,\lambda)$ is a $\bz_+$-graded $\lie g[t]$--module once we declare the grade of $w_{\nu,\lambda}$ to be zero. In the case when $\lambda=0$ it is known (see \cite{CPweyl}) that the  relations in \eqref{lnu} are a    consequence  of the relations in \eqref{localweyld}. In particular the module $M(\nu,0)$ is just the local Weyl module, which is usually denoted as $W_{\loc}(\nu)$. The local Weyl modules are known (see \cite{CPweyl}) to be  finite--dimensional and since $M(\nu,\lambda)$ is obviously a quotient of $W_{\loc}(\nu+\lambda)$, it follows that $M(\nu,\lambda)$ is also finite--dimensional. Moreover standard arguments show that \begin{equation*}\label{gl1}\dim\Hom_{\lie g}(V(\mu), M(\nu,\lambda))\ne 0\implies \nu+\lambda-\mu\in Q^+, \end{equation*}$$ \dim\Hom_{\lie g}(V(\nu+\lambda), M(\nu,\lambda))=1.$$
This discussion shows that the set
$\{\ch_{\gr}M(\mu,0): \mu\in P^+\}$ (resp.  the set $\{\ch_{\gr} M(0,\mu):\mu\in P^+\}$) is linearly independent and that the $\bz[q,q^{-1}]$--span  of the set contains
$\ch V(\lambda)$, $\lambda\in P^+$. It follows that  we can write, 
\begin{equation}\label{defgh} \ch_{\gr}M(\nu,\lambda)=\sum_{\mu\in P^+} g_{\nu,\lambda}^\mu(q)\ch_{\gr} M(\mu,0)= \sum_{\mu\in P^+} h_{\nu,\lambda}^\mu(q)\ch_{\gr} M(0,\mu), \end{equation} for some  $g_{\nu,\lambda}^\mu, h_{\nu,\lambda}^\mu\in\bc[q]$ satisfying, \begin{equation} \label{elempropgh} g_{\nu,\lambda}^{\nu+\lambda}=1= h_{\nu,\lambda}^{\nu+ \lambda},\ \ g_{\nu,\lambda}^\mu= h_{\nu,\lambda}^\mu= 0\ \ {\rm{if}}\ \  \lambda+\nu-\mu\notin Q^+.\end{equation} (In what follows we shall assume that $g_{\lambda,\mu}^\nu =0$ if one of $\lambda,\mu,\nu$ are not in $P^+$).
The linear independence of the  characters implies   that for all $\nu,\mu\in P^+$, \begin{equation}\label{inverse}\sum_{\mu'\in P^+}h_{\nu,0}^{\mu'}g_{0,\mu'}^\mu=\delta_{\nu,\mu}=\sum_{\mu'\in P^+}g_{0,\nu}^{\mu'}h_{\mu',0}^\mu.\end{equation}
It was shown in  \cite{CL}  that $W_{\loc}(\nu)$ (or equivalently $M(\nu,0)$) is graded isomorphic to a Demazure module occurring in a level one representation of the affine Lie algebra $\wh{\lie{sl}}_{n+1}$. Using \cite{San} it follows that \begin{equation*}\label{grweyl}\ch_{\gr} M(\lambda,0)= P_\lambda(z,q,0),\end{equation*} 
and hence we have 
\begin{equation}\label{grdem}\ch_{\gr} M(\nu,\lambda)= \sum_{{\mu\in P^+}} g_{\nu,\lambda}^\mu P_\mu(z,q,0).\end{equation} 
\begin{rem} The more usual definition of the Demazure modules in the affine case is by using affine  dominant integral weights and an element of the affine Weyl group. To make the connection with this definition, let $w_\circ$ be the longest element of the finite Weyl group and let $\Lambda_0$ be the affine dominant integral weight associated to the zero node of the affine Dynkin diagram. Choose $w$ in the affine Weyl group such that the element $w_\circ\nu+\Lambda_0=\Lambda$ is affine dominant integral. Then $M(\nu,0)$ corresponds to the Demazure module associated to the pair $(w,\Lambda)$.
\medskip

It is known  (see \cite{CSVW} and the references in that paper) 
that   the module $M(0,\lambda)$ is isomorphic to a Demazure module (denoted in the literature as $D(2,\lambda)$) occurring in a level two representation of the affine Lie algebra $\wh{\lie{sl}}_{n+1}$. In this case one makes the connection with the usual definition of Demazure modules by choosing an element $w$  of the affine Weyl group which makes $w_\circ\lambda+2\Lambda_0$ affine dominant integral.
\end{rem}

\subsection{} The following is the key result needed to prove the main theorem.
\begin{prop}\label{thmcrux} For $\lambda,\nu,\mu\in P^+$ with $(\nu,\lambda)$ admissible, we have \begin{gather} 
g_{\nu,\lambda}^\mu=q^{(\lambda+\nu-\mu,\nu)}g_{0,\lambda}^{\mu-\nu}\ \ \ \  h_{\nu,0}^\mu=p_\nu^\mu.
\end{gather}
\end{prop}
The proof of the first equality can be found in Section \ref{gnlm} while the proof of the second equality is in Section \ref{section4}. 

\subsection{Proof of Theorem \ref{main1a}} \label{strategy}
We now deduce Theorem \ref{main1a}.  Using the second equality in Proposition \ref{thmcrux} and equations \eqref{defg}, \eqref{defgh} we get $$\sum_{\mu\in P^+}p_\nu^\mu G_\mu(z,q)= P_\nu(z,q,0)= \ch_{\gr}M(\nu,0)=\sum_{\mu\in P^+}p_\nu^\mu \ch_{\gr}M(0,\mu).$$ Since $ M(\omega_i,0)\cong V(\omega_i)$ we have $G_{\omega_i}(z,q)= \ch_{\gr}M(0,\omega_i)$ and an induction  on the partial order on $P^+$ proves that $$G_\mu(z,q)=\ch_{\gr} M(0,\mu)\ \ {\rm{for \ all}} \ \mu\in P^+.$$  Using \eqref{defg} and \eqref{defgh} again we have $$\sum_{\mu\in P^+}a_\lambda^\mu P_\mu(z,q,0)=G_\lambda(z,q)=\sum_{\mu\in P^+}g_{0,\lambda}^\mu\ch_{\gr} M(\mu,0)=\sum_{\mu\in P^+}g_{0,\lambda}^\mu P_\mu(z,q,0).$$ The linear independence of the Macdonald polynomials implies   that \begin{equation}\label{a=g}a_\lambda^\mu= g_{0,\lambda}^\mu,\ \ {\rm{
for\  all}}\ \ \mu\in P^+.\end{equation}
Using the first equality in Proposition \ref{thmcrux} and equation  \eqref{gnl} we now have $$G_{\nu,\lambda}(z,q)=\ch_{\gr} M(\nu,\lambda)\ \ {\rm{if}}\ \ (\nu,\lambda)\in P^+\times P^+\ \ {\rm{admissible}}.$$  Since the modules $M(\nu,\lambda)$ are finite--dimensional and each graded piece of the module is a $\lie g$--module it follows that $G_{\nu,\lambda}(z,q)$ is Schur positive and hence  Theorem \ref{main1a} is proved. 
 
\subsection{} The proof of Proposition \ref{thmcrux} requires a deeper understanding of the modules $M(\nu,\lambda)$ and   we now state the key representation theoretic result of this paper.

 \begin{prop}\label{ses} 
\begin{enumerit}
 \item[(i)] Let $(\nu,\lambda)\in P^+\times P^+$ be an admissible pair. \\
 \begin{enumerit}
  \item[(a)]
If   $\nu(h_j)\ge 2$ for some $j\in[1,n]$ there exists an  exact sequence of $\lie g[t]$--modules, 
 $$0\to  \tau^*_{(\lambda_0+\nu,\alpha_j)-1}M(\nu-\alpha_j,\lambda)\to M(\nu,\lambda)\to M(\nu-2\omega_j,\lambda+2\omega_j)\to 0.$$
\item[(b)]
 If  $\nu\in P^+(1)$ with  $\max\nu<\min\lambda_1=m$ and  $p=\min(\lambda_1-\omega_m)>0$   we have an exact sequence of $\lie g[t]$--modules,$$0\to  \tau^*_{\frac12(\lambda, \alpha_{m,p})}M(\nu+\omega_{m-1},\lambda-\alpha_{m,p}-\omega_{m-1})\to M(\nu+\omega_m,\lambda-\omega_m)\to M(\nu,\lambda)\to 0.  $$
 \end{enumerit}
 \item[(ii)] If $  \lambda\in P^+(1)$ and $m\in[1,n]$ with $m<\min\lambda$, we have an exact  sequence of $\lie g[t]$- modules, $$ 0\to \tau^*_1 M(\omega_{m-1}, \lambda+\omega_{m+1})\to M(\omega_m, \lambda+\omega_m)\to M(0, \lambda+2\omega_m)\to 0.$$
 \end{enumerit} 
 \end{prop}

\medskip

\subsection{}\label{right}  From now on we prove Proposition \ref{ses}. As  a  first step we establish the  existence of the appropriate right exact sequences. Part (i) was proved in \cite{Wand} and we just indicate the main steps. For part (a) notice that since $(\nu,\lambda)$ is admissible it follows that if $\nu(h_j)\ge 2$ then either $\lambda_1=0$ or $\min\lambda_1>j+1$.
In this case the existence of the map $M(\nu,\lambda)\to M(\nu-2\omega_j,\lambda+2\omega_j)\to 0$ is trivial from the defining relations and also the fact that
 $(x_j^-\otimes t^{(\nu+\lambda_0,\alpha_j)-1})w_{\nu,\lambda}$ is in the kernel of the map. One then proves that this element  generates the kernel and that there exists a well-defined map $$M(\nu-\alpha_j,\lambda)\twoheadrightarrow \bu(\lie g[t])(x_j^-\otimes t^{(\nu+\lambda_0,\alpha_j)-1})w_{\nu,\lambda}.$$
The proof of part (b) is similar. This time one proves that there is  a canonical map $M(\nu+\omega_m,\lambda-\omega_m) \to M(\nu,\lambda)\to 0$ whose kernel is generated by $(x_{m,p}^-\otimes t^{\frac12(\lambda,\alpha_{m,p})})w_{\nu+\omega_m,\lambda-\omega_m}$. One then proves that  there exists a well-defined map $$M(\nu+\omega_{m-1},\lambda-\alpha_{m,p}-\omega_{m-1})\twoheadrightarrow \bu(\lie g[t])(x_{m,p}^-\otimes t^{\frac12(\lambda,\alpha_{m,p})})w_{\nu+\omega_m,\lambda-\omega_m}.$$

\medskip

The proof of part (ii) is similar but since this is not in the literature, we provide all the details for the reader's convenience. An inspection of the  defining relations shows the existence of  a surjective map of $\lie g[t]$-modules  $$M(\omega_m,\lambda+\omega_m)\to M(0,\lambda+2\omega_m)\to 0,$$ whose kernel is generated by the elements $$(*)\ \ \left(x_\alpha^-\otimes t^{\lceil(\lambda+\omega_m)(h_\alpha)/2\rceil}\right)w_{\omega_m,\lambda+\omega_m},\ \ \alpha\in R^+,\ \   \lambda(h_\alpha)\in 2\bz_+,\ \  \omega_m(h_\alpha)=1.$$ Taking $\alpha=\alpha_m$ we see that $(x_m^-\otimes t)w_{\omega_m,\lambda+\omega_m}$ is in the kernel and we claim that it generates the kernel. 
For the claim, let $\alpha\in R^+$ be as in $(*)$ and write  $$ \alpha=\beta+\alpha_m+\gamma,\ \ \  \beta\in \sum_{i=1}^{m-1}\bz_+\alpha_i,\ \  \gamma\in\sum_{i=m+1}^n\bz_+\alpha_i.$$  If $\beta,\gamma\in R^+$ then $(\lambda+\omega_m)(h_\beta)=0 $  and $(\lambda+\omega_m)(h_\gamma)=\lambda(h_\gamma)$; hence the  defining relations  of $M(\nu,\lambda)$ give
 $$(x_\beta^-\otimes 1)(x_\gamma^-\otimes t^{\lceil\lambda(h_\gamma)/2\rceil})(x_m^-\otimes t)w_{\omega_m,\lambda+\omega_m}=
(x_\alpha^-\otimes t^{\lceil(\lambda+\omega_m)(h_\alpha)/2\rceil})w_{\omega_m,\lambda+\omega_m}.$$ Otherwise either  $\beta=0$ or $\gamma=0$  and 
an identical argument  gives the claim in this case.

We now  show that the existence of the morphism  $$\tau_1^* M(\omega_{m-1},\lambda+\omega_{m+1})\to M(\omega_m,\lambda+\omega_m),\ \ w_{\omega_{m-1},\lambda+\omega_{m+1}}\mapsto (x_m^-\otimes t)w_{\omega_m,\lambda+\omega_m}.$$   The proof that $(x_m^-\otimes t)w_{\omega_m,\lambda+\omega_m}$  satisfies the relations in \eqref{localweyld} is elementary. It remains to show that $$(\dagger)\ \ (x^-_\alpha\otimes t^{\omega_{m-1}(h_\alpha)+\lceil{(\lambda+\omega_{m+1}})(h_\alpha)/2\rceil})(x_m^-\otimes t)w_{\omega_m,\lambda+\omega_m}=0,$$ and for this we
 consider several cases. If $\alpha=
\alpha_m$ the relations 
$$(h_m\otimes 1)(x_m^-\otimes t)w_{\omega_m,\lambda+\omega_m}=(\lambda+2\omega_m-\alpha_m)(h_m)(x_m^-\otimes t)w_{\omega_m,\lambda+\omega_m}=0$$
and$$ (x_m^+\otimes 1)(x_m^-\otimes t) w_{\omega_m,\lambda+\omega_m}=0$$
imply that  $$(x_m^-\otimes 1)(x_m^-\otimes t)w_{\omega_m,\lambda+\omega_m}=0$$   since $M(\omega_m,\lambda+\omega_m)$ is finite-dimensional $\lie g$-module.

If $\alpha$ is in the span of $\{\alpha_j: 1\le j<m-1\}$ or
 $\{\alpha_j: m+1<j\le n\}$ then $(\dagger)$ follows from the relations
$$(x_\alpha^-\otimes t^{\lceil\lambda(h_\alpha)/2\rceil})w_{\omega_m,\lambda+\omega_m}=0,\ \ [x_m^-,x_\alpha^-]=0.$$ 

Suppose that $\alpha+\alpha_m\in R^+$ and $\alpha$ is in the span of $\{\alpha_j: 1\le j\le m-1\}$ (resp.  $\{\alpha_j: m+1\le j\le n\}$) and we have $$[x_{\alpha}^-,x_m^-]=x_{\alpha+\alpha_m}^-,\ \ \omega_{m-1}(h_\alpha)=1,\ \ \omega_m(h_\alpha)=\omega_{m+1}(h_\alpha)=0,$$
$$ ({\rm{resp.}}\ [x_{\alpha}^-,x_m^-]=x_{\alpha+\alpha_m}^-,\ \ \omega_{m-1}(h_\alpha)=1,\ \ \ \omega_{m-1}(h_\alpha)=0= \omega_m(h_\alpha),\ \ \omega_{m+1}(h_\alpha)=1).$$ The defining relations of $M(\omega_m,\lambda+\omega_m)$ give $$(x_\alpha^-\otimes t^{1+\lceil\lambda(h_\alpha)/2\rceil}) w_{\omega_m,\lambda+\omega_m}=0= (x_{\alpha+\alpha_m}^-\otimes t^{(2+\lceil\lambda(h_\alpha)/2\rceil)})w_{\omega_m,\lambda+\omega_m}, $$  
and $(\dagger)$ follows in this case.  

Finally we must consider the case when $\alpha=\beta+\alpha_m+\gamma$ with $\omega_{m-1}(h_\beta)=1=\omega_{m+1}(h_\gamma)$. Then $x_\alpha^-$ and $x_m^-$ commute and since $$\omega_{m-1}(h_\alpha)+\lceil(\lambda+\omega_{m+1})(h_\alpha)/2\rceil=\omega_m(h_\alpha)+\lceil(\lambda+\omega_{m})(h_\alpha)/2\rceil$$ we see that $(\dagger)$ becomes the relation $$ (x_m^-\otimes t)(x_\alpha^-\otimes t^{\omega_m(h_\alpha)+\lceil(\lambda+\omega_{m})(h_\alpha)/2\rceil})w_{\omega_m,\lambda+\omega_m}= 0,$$ in $M(\omega_m,\lambda+\omega_m)$. This completes the proof of the existence of the right exact sequence in (ii).
\medskip

In particular, we have proved that
for $(\nu,\lambda)$ admissible  \begin{gather}\label{le}\dim M(\nu,\lambda)\le\dim M(\nu-2\omega_j,\lambda+2\omega_j)+\dim M(\nu-\alpha_j,\lambda),\ \nu(h_j)\ge 2,\end{gather}
and if   $\lambda=2\lambda_0+\lambda_1$ and $\nu\in P^+(1)$ with
with $m=\max\nu<\max\lambda_1=p$, then
\begin{gather}\label{le1}\dim M(\nu, \lambda)\le \dim M(\nu-\omega_m,\lambda+\omega_m)+\dim M(\nu-\omega_m+\omega_{m-1},\lambda-\omega_p+\omega_{p+1}). \end{gather}
If $\lambda\in P^+(1)$ and $m<\min\lambda$ then \begin{gather}\label{le2}\dim M(\omega_m,\lambda+\omega_m)\le \dim M(\omega_{m-1},\lambda+\omega_{m+1})+\dim M(0,\lambda+2\omega_m).\end{gather} 

\subsection{} To show  that the right exact sequences established  are exact it suffices to prove that equality holds in \eqref{le}, \eqref{le1} and \eqref{le2}. This requires some additional work.  We recall  some  results on the structure of $M(\nu,0)$ and $M(0,\lambda)$.
 Part (i) of the following proposition  was proved in  \cite{CL} while parts  (ii) and (iii) were proved in \cite{CSVW} and \cite{Wand} respectively.  Note that the $M(\nu,0)$ is denoted as $W_{\loc}(\nu)$ and $M(0,\lambda)$ is denoted as $D(2,\lambda)$ in those papers.  Part (iv) was proved in \cite{BC} and establishes the  exactness of the sequence   in  Proposition \ref{ses}(i)(b) when $\nu=0$.
\begin{prop}\label{dimbcsw} 
\begin{enumerit}
\item[(i)] For all
  $\nu_1,\nu_2\in P^+$ where $\nu=\nu_1+\nu_2$, we have  $$\dim M(\nu,0)=\dim M(\nu_1,0)\dim M(\nu_2,0).$$
\item[(ii)] Suppose that $\lambda,\mu\in P^+$ with $\lambda=2\lambda_0+\lambda_1$ are such that $\lambda_0-\mu\in P^+$. Then
$$\dim M(0,\lambda)=\dim M(0,\lambda-2\mu)\dim M(0,2\mu).$$

\item[(iii)] For $\lambda,\nu\in P^+$ we have $$\dim M(\nu,\lambda)\ge \dim M(\nu,0)\dim M(0,\lambda). $$
\item[(iv)] For $m\in[1,n]$ and $\lambda\in P^+(1)$ with $m<\min\lambda=p$, we have $$\dim M(\omega_m,\lambda)=\dim M(\omega_m,0)\dim M(0,\lambda)$$ and  there exists  a short exact sequence of $\lie g[t]$-modules
 $$0\to \tau^*_1 M(\omega_{m-1}, \lambda-\omega_p+\omega_{p+1})\to M(\omega_m,\lambda)\to M(0,\lambda+\omega_m)\to 0.$$
\end{enumerit}\hfill\qedsymbol
\end{prop}

 \subsection{}  \label{paradm} 

 Let $(\nu,\lambda)\in P^+\times P^+$ is an admissible pair. Then, 
 \begin{itemize}
\item if $ \nu(h_j)\ge 2$,
for some $j\in[1,n]$ the pairs $(\nu-2\omega_j,\lambda+2\omega_j)$ and $ (\nu-\alpha_j,\lambda)$ are admissible.
\item If  $\nu\in P^+(1)$ with  $\max\nu<\min\lambda_1=m$,  the  pair $(\nu+\omega_m,\lambda-\omega_m)$  is admissible. If in addition $p=\min(\lambda_1-\omega_m)>0$, then the pair $(\nu+\omega_{m-1},\lambda-\alpha_{m,p}-\omega_{m-1})$ is also admissible. 
\end{itemize}
We shall use this observation freely in the rest of the paper.
 Define a partial order on the set of admissible pairs as follows: $$(\nu',\lambda')\le
 (\nu,\lambda)\ {\rm{if}}\ \ \nu+\lambda-\nu'-\lambda'\in Q^+\setminus\{0\}\  {\rm{ or}}\ \ \nu'+\lambda'=\nu+\lambda\ {\rm{and}}\  \nu-\nu'\in P^+.$$
 We claim that the elements $\{(0,\omega_i): i\in[0,n]\}$ are  the  minimal elements with respect to this order. To see this,  suppose that $(\nu,\lambda)$ is   admissible. If  $\nu(h_j)\ge 2$ for some $j\in[1,n]$ the pair  $(\nu-2\omega_j,\lambda+2\omega_j)$ is admissible and less than $(\nu,\lambda)$. If $\nu\in P^+(1) $ and $m=\min\nu>0 $ then $(\nu-\omega_m,\lambda+\omega_m)$    is admissible and less than $(\nu,\lambda)$.  If $\nu=0$ and  $i\in[0,n]$ is such that  $\lambda-\omega_i\in Q^+\setminus\{0\}$  then $(0,\omega_i)<(0,\lambda)$ and hence the claim is proved.
 
 \medskip

\noindent The following result  shows that equality holds in \eqref{le}, \eqref{le1} and \eqref{le2} and completes the proof of Proposition \ref{ses}.
\begin{prop}\label{lb} \begin{enumerit}
 \item[(i)] 
For $(\nu,\lambda)$ admissible we have, $$\dim M(\nu,\lambda)=\dim M(\nu,0)\dim M(0,\lambda),$$   $$\dim M(\nu,\lambda)=\dim M(\nu-2\omega_j,\lambda+2\omega_j)+\dim M(\nu-\alpha_j,\lambda),\ \nu(h_j)\ge 2.$$
If  $\lambda=2\lambda_0+\lambda_1$ and $\nu\in P^+(1)$ with
with $m=\max\nu<\min\lambda_1=p$,
$$\dim M(\nu, \lambda)=\dim M(\nu-\omega_m,\lambda+\omega_m)+\dim M(\nu-\omega_m+\omega_{m-1},\lambda-\omega_p+\omega_{p+1}). $$
\item[(ii)] If $\lambda\in P^+(1)$ and $m<\min\lambda$ then \begin{gather*}\dim M(\omega_m,\lambda+\omega_m)=\dim M(\omega_m,0)\dim M(0,\lambda+\omega_m)\\= \dim M(\omega_{m-1},\lambda+\omega_{m+1})+\dim M(0,\lambda+2\omega_m).\end{gather*} 
\end{enumerit}
\end{prop}
 
\begin{pf} 
  The proof of (i)  proceeds by an induction on the partial order on admissible pairs.    Induction begins for the minimal elements, since  $$\dim M(0,\omega_i)=\dim M(\omega_i,0)=\dim V(\omega_i),\ \ i\in[0,n].$$
  Before proving the inductive step we note the following.
  Using Proposition \ref{dimbcsw}(i) and elementary representation theory, we get for all $j\in[1,n]$,
 \begin{eqnarray*}& \dim M(2\omega_j,0)&= (\dim V(\omega_j))^2=  \dim V(2\omega_j)+\dim V(\omega_{j-1})\dim V(\omega_{j+1})\\&&=\dim M(0,2\omega_j)+\dim M(2\omega_j-\alpha_j,0)\end{eqnarray*}

 Assume that the result holds for all admissible pairs $(\nu',\lambda')<(\nu,\lambda)$. Suppose that there exists $j\in[1,n]$ with $\nu(h_j)\ge 2$. Then  Proposition \ref{dimbcsw}(iii) and  \eqref{le} give \begin{gather*} \dim M(\nu,0)\dim M(0,\lambda)\le \dim M(\nu,\lambda)
 \le \dim M(\nu-2\omega_j,\lambda+2\omega_j)+\dim M(\nu-\alpha_j,\lambda).\end{gather*} Using the inductive hypothesis and part (ii) of Proposition \ref{dimbcsw} we get,
 \begin{gather*}\dim M(\nu-2\omega_j,\lambda+2\omega_j)+\dim M(\nu-\alpha_j,\lambda)\\ =\dim M(\nu-2\omega_j,0)\dim M(0,\lambda+2\omega_j)+\dim M(\nu-\alpha_j,0)\dim M(0,\lambda)\\=
\left (\dim M(\nu-2\omega_j,0)\dim M(0,2\omega_j)+\dim M(\nu-\alpha_j,0)\right)\dim M(0,\lambda)\\= \dim M(\nu,0)\dim M(0,\lambda), \end{gather*} as needed. 

Suppose that  $\nu\in P^+(1)$ with $m=\min\nu$.  If $\lambda=2\lambda_0$, the isomorphism in \eqref{iso1} gives $$ \dim M(\nu,2\lambda_0)=\dim M(\nu-\omega_m,2\lambda_0+\omega_m).$$ Since $(\nu-\omega_m,\omega_m+2\lambda_0)<(\nu,\lambda)$ the inductive hypothesis, equation \eqref{iso1} and Proposition \ref{dimbcsw}(ii) give \begin{gather*}\dim M(\nu-\omega_m,2\lambda_0+\omega_m)=\dim  M(\nu-\omega_m,0)\dim M(0,\omega_m+\lambda)\\=
 \dim  M(\nu-\omega_m,0)\dim M(\omega_m,0)\dim M(0,2\lambda_0)= \dim M(\nu,0)\dim M(0,2\lambda_0),\end{gather*}
 where the last equality follows from Proposition \ref{dimbcsw}(i). 
 
 If  $\nu\in P^+(1)$ with $\min\nu=m$ and   $p=\min\lambda_1>0$ 
  Proposition \ref{dimbcsw}(iii) and \eqref{le1} give \begin{gather*}\dim M(\nu,0)\dim M(0,\lambda)\le \dim M(\nu,\lambda)\le\\ \ \  \dim M(\nu-\omega_m,\lambda+\omega_m)+\dim M(\nu-\omega_m+\omega_{m-1},\lambda-\omega_p+\omega_{p+1}).\end{gather*}  
 The inductive hypothesis applies to  the modules on the right hand side of the preceding inequality 
and so we get
 \begin{gather*} \dim M(\nu-\omega_m,\lambda+\omega_m)= \dim M(\nu-\omega_m,0)\dim M(0,\lambda+\omega_m)\\= \dim M(\nu-\omega_m,0)\left(\dim M(\omega_m,0)\dim M(0,\lambda)- \dim M(\omega_{m-1},0)\dim M(0,\lambda-\omega_p+\omega_{p+1})\right)\\=\dim M(\nu,0)\dim M(0,\lambda)-\dim M(\nu-\omega_m+\omega_{m-1}, \lambda-\omega_p+\omega_{p+1})
 \end{gather*}where the second equality follows from Proposition \ref{dimbcsw}(ii) and (iv).  For the last equality we  use Proposition \ref{dimbcsw}(i) and the inductive hypothesis. This completes the proof of part (i) of this proposition.

   The proof of part (ii) is similar with \eqref{le2} showing  that  it is enough to prove 
$$\dim M(\omega_m,0)\dim M(0,\lambda+\omega_m)= \dim M(\omega_{m-1},\lambda+\omega_{m+1})+\dim M(0,\lambda+2\omega_m).$$
   We proceed by a downward induction on $m$. If $m=n$ then $\lambda=0$ by our assumptions and the result is known from our initial observations  since $M(\omega_m,\omega_m)\cong M(2\omega_m,0)$.
  For the inductive step,   Proposition \ref{dimbcsw}(iii) and \eqref{le2} give  \begin{gather*}\dim M(\omega_m,0)\dim M(0,\lambda+\omega_m)\le\dim M(\omega_m,\lambda+\omega_m)\\
  \le\dim M(\omega_{m-1},\lambda+\omega_{m+1})+\dim M(0,\lambda+2\omega_m)\end{gather*}
  Using Proposition \ref{dimbcsw}(ii), we have \begin{gather*}(*)\ \  \dim M(0,2\omega_m+\lambda)=(\dim M(2\omega_m,0)-\dim M(\omega_{m-1},0)\dim M(\omega_{m+1},0))\dim M(0,\lambda).\end{gather*}
   If we let $p=\min\lambda_1$ then  $(\omega_{m+1},\lambda)$ is admissible if $p>m+1$   and part (i) of this proposition applies; if $p=m+1$ the inductive hypothesis on $m$ applies to $(\omega_{m+1},\lambda)$. Hence $$\dim M(\omega_{m+1},0)\dim M(0,\lambda)=\dim M(0,\lambda+\omega_{m+1})+\dim M(\omega_m,\lambda-\omega_p+\omega_{p+1}).
   $$ Multiplying through by $\dim M(\omega_{m-1},0)$ and using the fact that $(\omega_{m-1},\lambda+\omega_{m+1})$ is admissible we have
  \begin{gather*}\dim M(\omega_{m-1},0) \dim M(\omega_{m+1},0)\dim M(0,\lambda)=\\\dim M(\omega_{m-1},\lambda+\omega_{m+1})+\dim M(\omega_{m-1},0)\dim M(\omega_m,\lambda-\omega_p+\omega_{p+1})
  \end{gather*}
   Substituting in $(*)$ and using the fact that $(\omega_m,\lambda-\omega_p+\omega_{p+1})$ is admissible gives
   \begin{gather*} \dim M(\omega_{m-1},\lambda+\omega_{m+1})+\dim M(0,\lambda+2\omega_m)\\=
\dim M(\omega_m,0)\left(\dim M(\omega_m,0)\dim M(0,\lambda)-\dim  M(\omega_{m-1},0)\dim M(0,\lambda-\omega_p+\omega_{p+1})\right)\\=\dim M(\omega_m,0)\dim M(0,\lambda+\omega_m)\end{gather*}
where the last equality follows since  $(\omega_m,\lambda)$ is admissible and we can  use part (i) of this proposition. \end{pf}

\section{The elements $g_{\nu,\lambda}^\mu$}\label{secfour}
We analyze recursive formulae for $g_{\nu,\lambda}^\mu$ which follow from Proposition  \ref{ses} and establish (see Section \ref{gnlm}) Proposition \ref{thmcrux} for the elements $g_{\nu,\lambda}^\mu$.
We also give closed formulae for the $g_{\nu,\lambda}^\mu$ and establish a technical result which will be needed in the study of the polynomials $h_{\nu,\lambda}^\mu$.

\subsection{} The next result is an immediate  consequence of Proposition \ref{ses}
 \begin{lem}\label{sescon} Assume that $(\nu,\lambda)$ is admissible. If $\nu(h_j)\ge 2$ for some $j\in[1,n]$ we have for $b\in\{g,h\}$ $$b^\mu_{\nu,\lambda}= b^\mu_{\nu-2\omega_j,\lambda+2\omega_j}+q^{(\lambda_0+\nu,\alpha_j)-1}b^\mu_{\nu-\alpha_j,\lambda},$$
 and if $\nu\in P^+(1)$ with $\max\nu<\min\lambda_1=m$ and $p=\min(\lambda_1-\omega_m)>0$, 
 $$b^\mu_{\nu+\omega_m,\lambda-\omega_m}=b^\mu_{\nu,\lambda} +q^{\frac12(\lambda,\alpha_{m,p})}b_{\nu+\omega_{m-1},\lambda-\alpha_{m,p}-\omega_{m-1`}}.$$
 Finally if $\lambda\in P^+(1)$ and $m<\min\lambda$ we have 
$$b_{\omega_m,\lambda+\omega_m}^\mu=b_{0,\lambda+2\omega_m}^\mu+qb_{\omega_{m-1},\lambda+\omega_{m+1}}^\mu.$$
 \end{lem}
 \begin{pf}
If $\nu(h_j)\ge 2$, then Proposition \ref{ses}(i)(a) gives an equality of graded characters \begin{equation*}
  \ch_{\gr} M(\nu,\lambda)=\ch_{\gr}M(\nu-2\omega_j,\lambda+2\omega_j)+ q^{(\nu+\lambda_0, \alpha_j)-1}\ch_{\gr}M(\nu-\alpha_j,\lambda).\end{equation*}  Using \eqref{defgh} and equating coefficients of $M(\mu,0)$ (resp. $M(0,\mu)$) on both sides gives the first assertion of the Lemma when $b=g$ (resp. $b=h$). The proof of the other assertions is identical.
 \end{pf}
  
\subsection{}\label{gnlm} We show by induction on $\Ht\lambda$ that if $(\nu,\lambda)$ is admissible, then $$g_{\nu,\lambda}^\mu=q^{(\lambda+\nu-\mu,\nu)}g_{0,\lambda}^{\mu-\nu}.$$ If $\Ht\lambda=0$ this follows since $g_{\nu,0}^\mu=\delta_{\nu,\mu}=g_{0,0}^{\mu-\nu}.$
For the inductive step we  consider various cases. 

Suppose first that $\lambda_1=0$ and let $j\in[1,n]$ be such that $(\lambda_0,\alpha_j)>0$.  Then $(\nu+2\omega_j,\lambda-2\omega_j)$ is admissible and we have, \begin{gather*}g_{\nu,\lambda}^\mu= g_{\nu+2\omega_j,\lambda-2\omega_j}^\mu-q^{(\lambda_0+\omega_j+\nu,\alpha_j)-1}g_{\nu+2\omega_j-\alpha_j,\lambda-2\omega_j}^{\mu}\\ =q^{(\lambda+\nu-\mu,\nu+2\omega_j)}g_{0,\lambda-2\omega_j}^{\mu-\nu-2\omega_j}-q^{(\lambda_0+\omega_j+\nu,\alpha_j)-1+(\lambda+\nu-\alpha_j-\mu,\nu+2\omega_j-\alpha_j)}g_{0,\lambda-2\omega_j}^{\mu-\nu-2\omega_j+\alpha_j}\\ =q^{(\lambda+\nu-\mu,\nu)}\left(g_{2\omega_j,\lambda-2\omega_j}^{\mu-\nu}- q^{(\lambda_0,\alpha_j)}g_{\omega_{j-1}+\omega_{j+1},\lambda-2\omega_j}^{\mu-\nu}\right)\\= q^{(\lambda+\nu-\mu,\nu)}g_{0,\lambda}^{\mu-\nu}.\end{gather*} Here the first equality follows from Lemma \ref{sescon}, the second and third follow from the inductive hypothesis applied to the admissible pairs $(\nu+2\omega_j,\lambda-2\omega_j)$, $(\nu+2\omega_j-\alpha_j,\lambda-2\omega_j)$ and  $(2\omega_j,\lambda-2\omega_j)$, $(2\omega_j-\alpha_j,\lambda-2\omega_j)$ respectively. The fourth equality is  a further application of Lemma \ref{sescon}.

If $\Ht\lambda_1=1$ and $m=\min\lambda_1$ then $(\nu+\omega_m,\lambda-\omega_m)$ is admissible. Since $\lambda-\omega_m= 2\lambda_0$ the  isomorphism in \eqref{iso1} and the inductive hypothesis  gives $$g_{\nu,\omega_m+2\lambda_0}^\mu=g_{\nu+\omega_m,2\lambda_0}^\mu=
 q^{(\lambda+\nu-\mu,\nu+\omega_m)}g_{0,2\lambda_0}^{\mu-\omega_m-\nu}
=q^{(\lambda+\nu-\mu,\nu)}g_{\omega_m, 2\lambda_0}^{\mu-\nu}= q^{(\lambda+\nu-\mu,\nu)}g_{0,\omega_m+ 2\lambda_0}^{\mu-\nu}.$$

For  $\Ht\lambda_1=r\ge 2$, let
 $m=\min\lambda_1$ and $p=\min(\lambda_1-\omega_m)$. Since $(\nu+\omega_m,\lambda-\omega_m)$
is admissible,  Lemma \ref{sescon} gives $$g_{\nu,\lambda}^\mu=g_{\nu+\omega_m,\lambda-\omega_m}^\mu-q^{(\nu+\frac12\lambda,\alpha_{m,p})}g_{\nu+\omega_{m-1},\lambda-\omega_m-\omega_p+\omega_{m+1}}^\mu.$$
The inductive hypothesis applies to both terms on the right hand side and so
\begin{eqnarray*} &g_{\nu,\lambda}^\mu&=q^{(\lambda+\nu-\mu,\nu+\omega_m)}g_{0,\lambda-\omega_m}^{\mu-\nu-\omega_m}-q^{(\nu+\frac12\lambda,\alpha_{m,p})+(\lambda+\nu-\alpha_{m,p}-\mu, \nu+\omega_{m-1})}g_{0,\lambda-\omega_m-\omega_p+\omega_{m+1}}^{\mu-\nu-\omega_{m-1}},\\
&&=q^{(\lambda+\nu-\mu,\nu)}\left(g_{\omega_m,\lambda-\omega_m}^{\mu-\nu}-q^{\frac12(\lambda,\alpha_{m,p})}g_{\omega_{m-1},\lambda-\omega_m-\omega_p+\omega_{p+1}}^{\mu-\nu}\right),\\ &&= g_{0,\lambda}^{\mu-\nu},
\end{eqnarray*}
where the second equality is again an application of the inductive hypothesis and the final equality follows from a further application of Lemma \ref{sescon}.

\subsection{}\label{recursion} We establish some further properties of the polynomials $g_{0,\lambda}^\mu$.
\begin{prop} For $\lambda\in P^+$ with $\lambda=2\lambda_0+\lambda_1$ and $m=\min\lambda_1,\ \ p=\min(\lambda-\omega_m)$,  we have,
\begin{gather}\label{gpropelem}  \sum_{\mu\in P}q^{\frac12(\mu,\mu)}g_{0,\lambda}^\mu=0,\\ \label{nonzero}  g_{0,\lambda}^\mu\ne 0\implies (\lambda-\mu,\omega_s)\le (\lambda,\alpha_s)\ \ {\rm{for \ all }}\ \ 1\le  s<m,
 \\\label{ht2ga} g_{0,\lambda}^\mu= g_{0,\lambda-\omega_m}^{\mu-\omega_m}-q(1-\delta_{p,0})g_{0,\lambda-\alpha_{m,p}-\omega_{m-1}}^{\mu-\omega_{m-1}},\ \ \lambda_0=0,\\ \label{2lambda0a} g_{0,\lambda}^\mu=g_{0,\lambda-2\omega_j}^{\mu-2\omega_j}-qg_{0,\lambda-2\omega_j}^{\mu-2\omega_j+\alpha_j},\ \ \lambda_0=\omega_j\ \ {\rm{and}}\ \  j+1<m.\end{gather}
\end{prop}
\begin{pf} If $m>0$ then the pair $(\omega_m,\lambda-\omega_m)$ is admissible and if $(\lambda_0,\alpha_j)> 0$  for some $j\in[1,n]$ then $(2\omega_j,\lambda-2\omega_j)$ is admissible. Hence,
Lemma \ref{sescon} and  the results of Section \ref{gnlm} show that 
 \begin{gather}
\label{ht2g} g_{0,\lambda}^\mu =q^{(\omega_m, \lambda-\mu)}g_{0,\lambda-\omega_m}^{\mu-\omega_m}-(1-\delta_{p,0})q^{\frac12(\lambda,\alpha_{m,p})+(\omega_{m-1},\lambda-\mu)}g^{\mu-\omega_{m-1}}_{ 0,\lambda-\alpha_{m,p}-\omega_{m-1}}\\ \label{2lambda0} g_{0,\lambda}^\mu=q^{(2\omega_j, \lambda-\mu)}g_{0,\lambda-2\omega_j}^{\mu-2\omega_j}-q^{(\lambda_0,\alpha_j)+(2\omega_j-\alpha_j,\lambda-\alpha_j-\mu)}g_{0,\lambda-2\omega_j}^{\mu-2\omega_j+\alpha_j}.\end{gather}
A straightforward induction on $\Ht\lambda$ proves \eqref{gpropelem} and \eqref{nonzero}. To establish \eqref{ht2ga}  it is enough to note that that \eqref{nonzero} applied to $(\lambda-\omega_m,\mu-\omega_m)$ and $(\lambda-\alpha_{m,p}-\omega_{m-1}, \mu-\omega_{m-1})$ gives $$g_{0,\lambda-\omega_m}^{\mu-\omega_m}\ne 0\implies (\lambda-\mu,\omega_m)=0,\ \ g^{\mu-\omega_{m-1}}_{ 0,\lambda-\alpha_{m,p}-\omega_{m-1}}\ne 0\implies (\lambda-\mu-\alpha_{m,p}, \omega_{m-1})=0.$$ The proof of \eqref{2lambda0a} is similar and we omit the details.

\end{pf}
 \subsection{} We shall need the following result in the next section.
 \begin{prop}\label{crucial} For $m\in[1,n]$ and $\lambda\in P^+(1)$ with $m<p=\min\lambda$ we have
\begin{equation}\label{crux1} g_{\omega_m, \lambda+\omega_m}^\mu=q^{(\lambda+ 2\omega_m-\mu,\  \omega_m)}g_{0,\lambda+\omega_m}^{\mu-\omega_m},\end{equation} or equivalently 
\begin{equation}\label{crux2} q^{(\lambda+2\omega_m-\mu,\ \omega_m)}g_{0,\lambda+\omega_m}^{\mu-\omega_m}=g_{0,\lambda+2\omega_m}^\mu+ q
g_{0,\lambda+\omega_{m+1}}^{\mu-\omega_{m-1}}.\end{equation}\end{prop}

\begin{pf} We first prove that \eqref{crux1} and  \eqref{crux2} are equivalent. The third equality in Lemma \ref{sescon}  gives $$g^\mu_{\omega_m,\lambda+\omega_m}= g^\mu_{0,\lambda+2\omega_m}+qg^\mu_{\omega_{m-1}, \lambda+\omega_{m+1}}.$$
Using Proposition \ref{thmcrux} on the admissible pair $(\omega_{m-1},\lambda+\omega_{m+1})$ shows that the previous equation can be be rewritten as, $$g_{\omega_m,\lambda+\omega_m}^\mu= g_{0,\lambda+2\omega_m}^\mu+ q^{1+(\omega_{m-1},\ \lambda+2\omega_m-\alpha_m-\mu)}g_{0,\lambda+\omega_{m+1}}^{\mu-\omega_{m-1}}.$$
Using \eqref{nonzero} we see that  the second term is non-zero only if  $(\omega_{m-1},\lambda+2\omega_m-\alpha_m-\mu)=0$ and hence we get,
$$g_{\omega_m,\lambda+\omega_m}^\mu= g_{0,\lambda+2\omega_m}^\mu+ qg_{0,\lambda+\omega_{m+1}}^{\mu-\omega_{m-1}}.$$ The equivalence of the two assertions is immediate.

We prove \eqref{crux2} by a downward induction on $m$. If $m=n$ then $\lambda=0$ and  $$g_{0,\ \omega_m}^{\mu-\omega_m}=\delta_{\mu-\omega_m,\omega_m},\ \ g_{0,0}^{\mu-\omega_{m-1}}=\delta_{\mu-\omega_{m-1},0}. $$ Hence when $m=n$ equation 
\eqref{crux2} follows if $$\delta_{\mu,2\omega_n}= g_{0,2\omega_n}^\mu+q \delta_{\mu,\omega_{n-1}}.$$ But this is precisely \eqref{2lambda0} applied to $2\omega_n$. For the inductive step we use
  \eqref{ht2ga} applied to the admissible pair $(\omega_m,\lambda)$ and either \eqref{ht2ga} or the inductive hypothesis to  $(\omega_{m+1},\lambda)$ to get
\begin{eqnarray*}
&g_{0,\lambda+\omega_m}^{\mu-\omega_m}&=  g_{0,\lambda}^{\mu-2\omega_m}- qg_{0,\lambda-\omega_p+\omega_{p+1}}^{\mu-\omega_m-\omega_{m-1}}, \\ 
&g_{0,\lambda+\omega_{m+1}}^{\mu-\omega_{m-1}}&= g_{0,\lambda}^{\mu-2\omega_m+\alpha_m}-qg_{0,\lambda-\omega_p+\omega_{p+1}}^{\mu-\omega_{m-1}-\omega_m},\end{eqnarray*} and hence we must prove that $$g_{0,\lambda+2\omega_m}^\mu=\left(q^{(\lambda+2\omega_m-\mu,\ \omega_m)}g_{0,\lambda}^{\mu-2\omega_m}-qg_{0,\lambda}^{\mu-2\omega_m+\alpha_m} \right) -q(q^{(\lambda+2\omega_m-\mu,\ \omega_m)}-q)g_{0,\lambda-\omega_p+\omega_{p+1}}^{\mu-\omega_{m-1}-\omega_m}. $$
Using \eqref{nonzero} we see that $$g_{0,\lambda-\omega_p+\omega_{p+1}}^{\mu-\omega_{m-1}-\omega_m}\ne 0\implies  (\lambda-\mu+2\omega_m-\alpha_{m,p}, \omega_m)=0\implies (\lambda+2\omega_m-\mu,\omega_m)=1,$$
$$g_{0,\lambda}^{\mu-2\omega_m}\ne 0\implies  (\lambda+2\omega_m-\mu,\omega_m)=0,$$ and hence it suffices to prove that $$g_{0,\lambda+2\omega_m}^\mu=g_{0,\lambda}^{\mu-2\omega_m}-qg_{0,\lambda}^{\mu-2\omega_m+\alpha_m} .$$
If $p>m+1$  the pair  $(2\omega_m,\lambda)$ is admissible and the preceding equation is precisely
 \eqref{2lambda0}. 
Suppose that  $p=m+1$ and $\lambda=\omega_{m+1}$ we have $$q^{(\omega_{m+1},\omega_{m+1}+2\omega_m-\mu) }g_{0,2\omega_m}^{\mu-\omega_{m+1}}=g_{\omega_{m+1}, 2\omega_m}^{\mu}= g_{0,\omega_{m+1}+2\omega_m}^\mu.$$ 
Hence we must prove that    $$\delta_{\omega_{m+1}, \mu-2\omega_m} =\delta_{\mu- \omega_{m+1}, 2\omega_m}-q\delta_{\mu-\omega_{m+1}, 2\omega_m-\alpha_m}+ q\delta_{\omega_{m+1}, \mu-2\omega_m+\alpha_m}, $$  and this is clear. Finally suppose that $p=m+1$ and 
 $r=\min(\lambda-\omega_p)>0$. Using the fact that $(\omega_{m+1},\lambda-\omega_{m+1})$, $(2\omega_m,\lambda-\omega_{m+1})$ and $(2\omega_m,\lambda-\omega_{m+1}-\omega_r+\omega_{r+1})$ are admissible, and using Lemma \ref{sec1} and Section \ref{gnlm} gives the following equalities: \begin{gather*}g_{0,\lambda}^{\mu-2\omega_m}=  g_{0,\lambda-\omega_{m+1}}^{\mu-2\omega_{m}-\omega_{m+1}}
- q g_{0,\lambda-\omega_{m+1}-\omega_r+\omega_{r+1}}^{\mu-3\omega_m}=g_{2\omega_m,\lambda-\omega_{m+1}}^{\mu-\omega_{m+1}}-q g_{2\omega_m,\lambda-\omega_m-\alpha_{m+1,r}}^{\mu-\omega_m},\\
= \left(g_{0,\lambda+2\omega_m-\omega_{m+1}}^{\mu-\omega_{m+1}}+ qg_{0,\lambda-\omega_{m+1}}^{\mu-2\omega_{m+1}-\omega_{m-1}}\right)-\left(g_{0,\lambda+\omega_m-\alpha_{m+1,r}}^{\mu-\omega_m}+ q g_{0,\lambda-\omega_m-\alpha_{m+1,r}}^{\mu-\omega_m-\omega_{m-1}-\omega_{m+1}}\right)\\
=\left(g_{0,\lambda+2\omega_m-\omega_{m+1}}^{\mu-\omega_{m+1}}-qg_{0,\lambda+\omega_m-\alpha_{m+1,r}}^{\mu-\omega_m}\right)+q\left(g_{\omega_{m+1},\lambda-\omega_{m+1}}^{\mu-\omega_{m+1}-\omega_{m-1}}-qg_{\omega_m,\lambda-\omega_m-\alpha_{m+1,r}}^{\mu-\omega_{m-1}-\omega_{m+1}}\right)
\end{gather*}
A further application of \eqref{2lambda0a} to the first term and \eqref{ht2ga} to the second term show that the right hand side of the last equality is precisely the right hand side of \eqref{crux2} and the inductive step is established and the proof of the proposition is complete.


\end{pf}

 \subsection{A closed form for $g_{0,\lambda}^\mu$} We now give a closed form for the elements $g_{0,\lambda}^\mu$. Notice that these elements are defined completely by the recursive formulae in \eqref{ht2g} and \eqref{2lambda0} with the initial condition $g_{0,0}^\mu=\delta_{\mu,0}$.  Let $\rho=\frac12\sum_{\alpha\in R^+}\alpha$  and for $\eta=\sum_{i=1}^ns_i\alpha_i\in Q^+$ set $\eta^\vee=\sum_{i=1}^ns_i\omega_i.$
 A simple checking using the $q$-binomial identity $$\qbinom{m}{s}=q^{s}\qbinom{m-1}{s}+\qbinom{m-1}{s-1},\ \ m,s\in\bz_+$$  now shows that equation \eqref{2lambda0} forces, 
 \begin{eqnarray}\label{qbinom1} &g^\mu_{0,2\lambda_0}&=(-1)^{(2\lambda_0-\mu,\rho)} q^{\frac{1}{2}(2\lambda_0-\mu,\ \mu+\rho+(2\lambda_0-\mu)^\vee)}\prod_{i=1}^n\qbinom{(\lambda_0,\alpha_i)}{(2\lambda_0-\mu, \omega_i)}.\end{eqnarray}
 Equation \eqref{ht2ga} now gives a closed form for $g_{0,\lambda}^\mu$ when $\Ht\lambda_1=1$. To give a closed form for $g_{0,\lambda}^\mu$ when $\Ht\lambda_1\ge 2$ needs some additional work.

 For $s\ge 0$ and $\lambda\in P^+(1)$ define subsets $\Sigma_s(\lambda)=\Sigma_s^0(\lambda)\cup\Sigma_s^1(\lambda)$ of $P^+$ inductively by:
\begin{gather*} \Sigma_0^0(\lambda)=\lambda,\ \ \Sigma_0^1(\lambda)=\emptyset,  \ \ \Sigma_s^0(\omega_m)=\Sigma_s^1(\omega_m)=\emptyset,\ \ m\in[0,n], \ s>0, \end{gather*} 
and for $s>0$ and  $\Ht\lambda\ge 2$ with $m=\min\lambda$, $p=\min(\lambda-\omega_m)$, we take  
\begin{eqnarray*}
&\Sigma_s^0(\lambda)&=\omega_m+\Sigma_s(\lambda-\omega_m),\\
&\Sigma^1_s(\lambda)&=\Sigma^0_{s-1}(\lambda-\alpha_{m,p}),\ \ \lambda(h_{p+1})=0\ \ 
\end{eqnarray*}
  and if $\lambda(h_{p+1})=1$ then $\Sigma^1_1(\lambda)=\{\lambda-\alpha_{m,p}\}$ and if  $s\ge 2$  $$\Sigma_s^1(\lambda)= (2\omega_{p+1}+\Sigma^0_{s-1}(\lambda-2\omega_{p+1}-\alpha_{m,p}))\cup (2\omega_{p+1}-\alpha_{p+1}+\Sigma^0_{s-2}(\lambda-2\omega_{p+1}-\alpha_{m,p})).$$
 Here we understand  that if  $\Sigma_s^r(\lambda)=\emptyset$ for $r\in\{0,1\}$ then so is  the set $\nu + \Sigma_s^r(\lambda) $.

\vskip12pt
We give a few examples of the sets $\Sigma_s(\lambda)$.  For $m<p$, 
$$\Sigma_s^0(\omega_m+\omega_p)=\begin{cases}\{\omega_m+\omega_p\},\ \ s=0,\\ \emptyset,\ \  s>0,\ \end{cases} \ \ \Sigma^1_s(
 \omega_m+\omega_p)=\begin{cases} \{\omega_{m-1}+\omega_{p+1}\},\ \ s=1,\\
 \emptyset,\  \ s\ne 1.\end{cases}
 $$
 For $m<p<r$,
 \begin{eqnarray*}&\Sigma_s^0(\omega_m+\omega_p+\omega_r)&=\begin{cases}\{\omega_m+\omega_p+\omega_r,\}\ \ s=0,\\
 \{\omega_{m}+\omega_{p-1}+\omega_{r+1}\},\ \ s=1,\\ \emptyset,\ \  {\rm{otherwise}},\ \end{cases}\\ \\ &\Sigma^1_s(\omega_m+\omega_p+\omega_r)&=\begin{cases} \{\omega_{m-1}+\omega_{p+1}+\omega_r\},\ \ \ s=1,\\
 \{\omega_{m-1}+\omega_p+\omega_{r+1}\},\ \ s=2,\\
 \emptyset\  \ {\rm{otherwise}}.\end{cases}
 \end{eqnarray*}
 Our final example is $\lambda=\omega_m+\omega_p+\omega_{\ell}+\omega_r$ with $m<p<\ell<r$. This is the first time we see the dependence on $\lambda(h_{p+1})$ and the non-empty sets are listed below:
 \begin{eqnarray*}&\Sigma^0_s(\lambda)&=\begin{cases}\{\omega_m+\omega_p+\omega_\ell+\omega_r\},\ \ s=0,\\
 \{\omega_m+\omega_{p-1}+\omega_{\ell+1}+\omega_r,\ \ \omega_m+\omega_p+\omega_{\ell-1}+\omega_{r+1}\},\ \ s=1,\\
\{ \omega_m+\omega_{p-1}+\omega_\ell+\omega_{r+1}\},\ \ s=2,\\
 \end{cases}\\ \\
 &\Sigma^1_s(\lambda)&=\begin{cases}\{\omega_{m-1}+\omega_{p+1}+\omega_\ell+\omega_r\},\ \  s=1,\\
 \{\omega_{m-1}+\omega_p+\omega_{\ell+1}+\omega_r,\   \omega_{m-1}+\omega_{p+1}+\omega_{\ell-1}+\omega_{r+1}\}\ \ s=2, \ell\ne p+1\\
 \{\omega_{m-1}+\omega_p+\omega_{\ell+1}+\omega_r\}\,\ \ s=2\ \ {\rm{and}}\ \ \ell= p+1,\\
 \{\omega_{m-1}+\omega_p+\omega_\ell+\omega_{r+1}\},\ \ s=3\ \ {\rm{and}}\ \ \ell\ne p+1.\\ 
 \end{cases}
 \end{eqnarray*}
 
 \begin{lem}\label{sigmaset} Let $\lambda\in P^+(1)$ and $\mu\in\Sigma_s(\lambda)$ for some $s
 \ge 0$. Then $$\mu(h_k)=0,\ \  k<\min\lambda-1,\ \ \mu(h_{\min\lambda-1})\le 1  \  {\rm{and}}\ \  \lambda-\mu\in\sum_{k\ge \min\lambda}\bz_+\alpha_k.$$ In particular the sets $\Sigma_s^r(\lambda)$ and $\Sigma_{s'}^{r'}(\lambda)$ are disjoint unless $(s,r)=(s',r')$.
 \end{lem}
 \begin{pf}
 The proof of the displayed statements is immediate from the definition of $\Sigma_s(\lambda)$ and a straightforward induction on $\Ht\lambda$.  The same induction also shows that the sets $\Sigma_s^r(\lambda)$ and $\Sigma_{s'}^{r
 }(\lambda)$ are disjoint for $r\in\{0,1\}$. Moreover if $m=\min\lambda$, $\mu\in\Sigma_s^0(\lambda)$, say $\mu=\omega_m+\nu$, then $\mu(h_{m-1})=0$ and hence $\mu\notin \Sigma_{s'}^{r'}(\lambda)$.
 \end{pf}

  For $\lambda,\mu\in P^+$, with $\lambda=2\lambda_0+\lambda_1$ and $\lambda_1\ne 0$ we have
 \begin{eqnarray}
 \label{geng} &g_{0,\lambda}^\mu&=q^{\frac12(\lambda_1, \lambda)}\sum_{s\ge 0}(-1)^s\sum_{\nu\in\Sigma_s(\lambda_1)}
  q^{\frac12(2\lambda_0-2\mu+\nu,\nu)}g_{0,2\lambda_0}^{\mu-\nu},\ \ 
 \end{eqnarray}
 where we understand that the second summation is zero if $\Sigma_s(\lambda_1)=\emptyset$. 
 The proof is a tedious calculation whch amounts to establishing that the right hand side of \eqref{geng} satisfies equation \eqref{ht2g}.

  \section{The polynomials $h_{\nu,0}^\mu(q)$}\label{section4}
 \subsection{}   Recall from Section \ref{sec1} that for $\nu,\mu\in P^+$ with $\mu=2\mu_0+\mu_1$, we defined \begin{equation*}\label{qbinom}p_\nu^\mu=  q^{\frac12(\nu+\mu_1,\  \nu-\mu)}\displaystyle\prod_{j=1}^n\qbinom{(\nu-\mu,\ \omega_j) +( \mu_0,\alpha_j)}{(\nu-\mu,\omega_j)}\ \ \end{equation*}
 and $p_\nu^\mu=0$ if $\nu$ or $\mu$ is not in $P^+$. It is easy to check that
 $$p_\nu^\mu= q^{(\omega_j,\nu-\mu)}p_{\nu-2\omega_j}^{\mu-2\omega_j}+q^{(\nu,\alpha_j)-1}p_{\nu-\alpha_j}^\mu  \ \ \  \  j\in[1,n],\ \ \nu(h_j)\ge 2,$$ $$ p_{\nu+\omega_m}^{\mu+\omega_m}=q^{(\nu-\mu,\omega_m)}p_\nu^\mu,\ \ {\rm{if}}\ \ \mu(h_m)\in2\bz_+.\ \ $$ 
In the rest of this section we will prove
that  \begin{equation}\label{final} h_{\nu,0}^\mu =p_\nu^\mu  \ \ {\rm{for\ all}} \ \ \nu,\mu\in P^+.\end{equation}
which is  the second assertion of  Proposition \ref{crux}.

 \subsection{} The next Lemma proves  \eqref{final} when $\Ht\mu\le 1$.
\begin{lem}\label{hinit} Let  $\nu\in P^+$ and  $r\in[0,n]$ be the  unique integer satisfying  $\nu-\omega_r\in Q^+$. Then \begin{gather*}
h_{\nu,0}^{\omega_k}   = q^{\frac12(\nu+\omega_r,\nu-\omega_r)}\delta_{k,r},\ \ k\in[0,n].\end{gather*}
\end{lem}
\begin{pf} The lemma is proved by induction on $\nu$ with respect to the partial order on $P^+$. If $\nu=\omega_k$ for some $k\in[0,n]$  then $\omega_k-\omega_r\in Q^+$ iff $k=r$ and the  result follows since $h_{\nu,0}^\mu=0$ if $\nu-\mu\notin  Q^+$ (see \eqref{elempropgh}).
Assume the result for all $\nu'\in P^+$ with $\omega_r\le \nu'<\nu$, in particular we have  $\Ht\nu\ge 2$. Using \eqref{inverse}  and the inductive hypothesis gives,$$ (*)\ \  \  0 = \sum_{{\mu\in P^+}} g_{0,\nu}^{\mu} h_{\mu,0}^{ \omega_k}= h_{\nu,0}^{\omega_k} +\sum_{\omega_k< \mu} g_{0,\nu}^{\mu}h_{\mu,0}^{\omega_k}= h_{\nu,0}^{\omega_k} +\delta_{k,r}\sum_{{\mu<\nu}}g_{0,\nu}^{\mu} q^{\frac12(\mu+\omega_r,\  \mu-\omega_r)}.$$  Equation \eqref{gpropelem} gives $$q^{\frac12(\nu+\omega_r,\nu-\omega_r)}+\sum_{{\mu<\nu}}q^{\frac12(\mu+\omega_r,\mu-\omega_r)}g_{0,\nu}^\mu=0,$$ and substituting in $(*)$ proves the Lemma.  \end{pf}

\subsection{} We need several properties of the elements $h_{\nu,\lambda}^\mu$. \begin{lem}\label{*} Suppose that $(\nu,\lambda)\in P^+\times P^+$ is admissible or that  $(\nu,\lambda)=(\omega_m,\lambda)$ with $\min\lambda=m$. Then
 \begin{equation*}\label{crux} \ \   h_{\nu,\lambda}^\mu=\sum_{\mu'\in P^+} q^{(\lambda+\nu-\mu',\nu)}g_{0,\lambda}^{\mu'-\nu}h_{\mu',0}^{\mu}.\ \ 
\end{equation*}
\end{lem}
\begin{pf} Recall from Section \ref{mnl} that $$\ch_{\gr}M(\nu,\lambda)=\sum_{\mu'\in P^+}g_{\nu,\lambda}^{\mu'}\ch_{\gr} M(\mu',0)=\sum_{\mu',\mu\in P^+}g_{\nu,\lambda}^{\mu'}h_{\mu',0}^\mu\ch_{\gr} M(0,\mu).$$  If  $(\nu,\lambda)$ is admissible  (resp. $(\nu,\lambda)=(\omega_m,\lambda+\omega_m)$) then using Proposition \ref{gnlm} (resp. Proposition \ref{crucial}) gives $$\sum_{\mu\in P^+}h_{\nu,\lambda}^\mu\ch_{\gr} M(0,\mu)=\ch_{\gr}M(\nu,\lambda)=\sum_{\mu',\mu\in P^+} q^{(\lambda+\nu-\mu', \nu)} g_{0,\lambda}^{\mu'-\nu}h_{\mu',0}^\mu\ch_{\gr} M(0,\mu).$$
The lemma follows by equating coefficients of $\ch_{\gr} M(0,\mu)$.\end{pf}

\subsection{} Suppose that $(\nu,\lambda)$ is admissible.   Recall from Lemma \ref{sescon} that, \begin{equation}\label{rec1ah}
  h_{\nu,\lambda}^\mu=h_{\nu-2\omega_j,\lambda+2\omega_j}^\mu+ q^{(\nu+\lambda_0,\alpha_j)-1}h_{\nu-\alpha_j,\lambda}^\mu,\ \ {\rm{if}}\ \ \nu(h_j)\ge 2,\end{equation}
  and if $\nu\in P^+(1)$ with $\max\nu=m<\min\lambda_1=p$ then
  \begin{equation}\label{rec1bh}   h_{\nu,\lambda}^\mu=h_{\nu-\omega_m, \lambda+\omega_m}^\mu+ q^{\frac12( \lambda,\alpha_{m,p})}h_{\nu-\omega_m+\omega_{m-1},\lambda-\omega_p+\omega_{p+1}}^\mu.\end{equation}

 \begin{lem} For $(\nu,\lambda)$ admissible and $k\in[1,n]$, we have $$h_{\nu,\lambda+2\omega_k}^{\mu+2\omega_k}=q^{(\lambda+\nu-\mu,\omega_k)}h_{\nu,\lambda}^\mu,$$ and hence for $j\in[1,n]$ with $ \nu(h_j)\ge 2$, 
 \begin{equation}\label{peeltwo}h_{\nu,\lambda}^\mu= q^{(\lambda+\nu-\mu,\omega_j)}h_{\nu-2\omega_j,\lambda}^{\mu-2\omega_j}+ q^{(\nu+\lambda_0,\alpha_j)-1}h_{\nu-\alpha_j,\lambda}^\mu.\end{equation}  
 \end{lem}
 \begin{pf} The proof is by induction on the partial order on admissible pairs (see Section \ref{paradm}). Induction obviously begins for the minimal elements $(0,\omega_i)$, $i\in[0,n]$ since $h_{0,\lambda}^\mu=\delta_{\lambda,\mu}$ for all $\lambda,\mu\in P^+$. Applying the inductive hypothesis to the right hand side of \eqref{rec1ah} we get if $\nu(h_j)\ge 2$ for some $j\in[1,n]$, $$h_{\nu,\lambda}^\mu= q^{-(\lambda+\nu-\mu,\omega_k)}\left(h_{\nu-2\omega_j,\lambda+2\omega_j+2\omega_k}^{\mu+2\omega_k}+ q^{(\nu+\lambda_0+\omega_k, \alpha_j)-1}h_{\nu-\alpha_j,\lambda+2\omega_k}^{\mu+2\omega_k}\right)= q^{-(\lambda+\nu-\mu,\omega_k)}h_{\nu,\lambda+2\omega_k}^{\mu+2\omega_k}.$$
 where the second equality follows by using \eqref{rec1ah} again  on the admissible  pair $(\nu,\lambda+2\omega_k)$.
 If $\nu\in P^+(1)$ then the inductive hypothesis applied to the right hand side of \eqref{rec1bh} gives \begin{eqnarray*}&h_{\nu,\lambda}^\mu&= q^{-(\lambda+\nu-\mu,\omega_k)}\left(h_{\nu-\omega_m,\lambda+\omega_m+2\omega_k}^{\mu+2\omega_k} +q^{( \lambda_0+\omega_k,\alpha_{m,p})+1}h_{\nu-\omega_m+\omega_{m-1},\lambda+2\omega_k-\omega_p+\omega_{p+1}}^{\mu+2\omega_k}\right)\\&&=q^{-(\lambda+\nu-\mu,\omega_k)}h_{\nu,\lambda+2\omega_k}^{\mu+2\omega_k},\end{eqnarray*} where the last equality is a further application of \eqref{rec1bh} to the admissible pair $(\nu, \lambda+2\omega_k)$. Equation \eqref{peeltwo} is now immediate  from \eqref{rec1ah} and  the proof of the Lemma is complete. \end{pf}

\subsection{} 

\begin{lem}\label{hmueven}
Let $m\in[1,n]$ and  $\lambda=2\lambda_0+\lambda_1\in P^+$ with $m<\min\lambda_1$ if $\lambda_0\ne 0$ and $m\le \min\lambda_1$ if $\lambda_0=0$.    Then
 $$h_{\omega_m,\lambda}^{\mu}=0\ \ {\rm{for}}\ \ \mu\in P^+ \ {\rm{with}}\  \mu(h_m)\in 2\bz_++1\ {\rm{and}}\  \ \lambda+\omega_m\ne\ \mu.$$ 
 \end{lem}

\begin{pf}
 If $m<\min\lambda_1$ we claim that   a stronger statement is true; namely$$ h_{\omega_m,\lambda}^{\mu}\ne 0\implies \mu=\lambda+\omega_m\ \ {\rm{or}}\ \ \mu(h_s)\in 2\bz_+\ \ {\rm{for\ all}}\  s\in[m,p].$$ 
We prove the claim  by an induction on $m$ with induction beginning  at $m=0$ since $h_{0,\lambda}^\mu=\delta_{\lambda,\mu}$.
  If $m\ge 1$ then using \eqref{rec1bh} we have, $$h_{\omega_m,\lambda}^\mu=0\implies \mu=\lambda+\omega_m\ \ {\rm{or}}\ \  h_{\omega_{m-1},\lambda-\omega_p+\omega_{p+1}}^\mu\ne 0.$$ The inductive hypothesis  applies to  the pair $(\omega_{m-1} ,\lambda-\omega_p+\omega_{p+1}) $ and so $$h_{\omega_{m-1},\lambda-\omega_p+\omega_{p+1}}^\mu\ne 0\implies\mu=\lambda+\omega_m-\alpha_{m,p}\ \ \ {\rm{or}}\ \ \mu(h_s)\in 2\bz_+\ {\rm{for\ all}}\ \ s\in[m-1, p+1].$$ Since $(\lambda_1+\omega_m-\alpha_{m,p})(h_s)=0$ for all $s\in[m,p]$ the inductive step follows and the claim is proved.
 \medskip

 It remains to prove the lemma when  $\lambda\in P^+(1)$ and  $m=\min\lambda_1$. Then Lemma \ref{sescon} gives  $$h_{\omega_m,\lambda}^{\mu}\ne 0\implies \mu=\lambda+\omega_m \ {\rm{or}}\ \  h^{\mu}_{\omega_{m-1},\lambda-\omega_m+\omega_{m+1}}\ne 0.$$
 If $\lambda-\omega_m=\omega_{m+1}$ then  the isomorphism in \eqref{iso1} gives $$h_{\omega_{m-1}, 2\omega_{m+1}}^\mu\ne 0\iff \mu=\omega_{m-1}+2\omega_{m+1},$$ and the result follows in this case. 
 Otherwise $\lambda-\omega_m\ne \omega_{m+1}$ and the result has been proved for $(\omega_{m-1},\lambda-\omega_m+\omega_{m+1})$.Hence  $$ h^{\mu}_{\omega_{m-1},\lambda-\omega_m+\omega_{m+1}}\ne 0\implies \mu(h_m)\in 2\bz_+\ \ {\rm{or}}\ \ \mu=\lambda-\omega_m+\omega_{m+1}+\omega_{m-1}.$$ Since $(\lambda-\omega_m+\omega_{m+1}+\omega_{m-1})(h_m)=0$ the proof of the Lemma is complete.

 \end{pf}

\subsection{Proof of equation \eqref{qbinom}}
The equality obviously holds if $\nu-\mu\notin Q^+$. If   $\nu-\mu=\sum_{i=1}^ns_i\alpha_i\in Q^+$ we set $\Ht_r(\nu-\mu)=\sum_{i=1}^ns_i$
and proceed by induction on $\Ht_r(\nu-\mu)$.  The induction begins when $\Ht_r(\nu-\mu)=0$ since then $\nu=\mu$ and  $h_{\nu,0}^\nu =1$. Assume the result holds for all pairs $(\nu,\mu)$ with $\Ht_r(\nu-\mu)<N$.
\medskip{}

We prove the result for $\Ht_r(\nu-\mu)=N$ by a  further induction on 
$\Ht\mu$.   Lemma \ref{hinit} shows that this second  induction begins when $\Ht\mu\le 1$. Assume the result holds for all $\nu$ with $\Ht_r(\nu-\mu)\le N$ and 
$1\le \Ht\mu<s$. 
If $\nu(h_j)\ge 2$ for some $j\in[1,n]$ then the inductive hypothesis on $\Ht\mu$ (resp. $\Ht_r(\nu-\mu)$) applies to the first (resp. second term) on the right hand  sides of \eqref{peeltwo} with ($\lambda=0$). Hence $$h_{\nu,0}^\mu=q^{(\omega_j,\nu-\mu)}p_{\nu-2\omega_j}^{\mu-2\omega_j}+q^{(\nu,\alpha_j)-1}p_{\nu-\alpha_j}^\mu=p_\nu^\mu.$$ In particular  the  inductive step has been proved for $\Ht \mu\in\{s,s+1\}$ provided that there exists $j\in[1,n]$ with $\nu(h_j)\ge 2$. 

If $\nu\in P^+(1)$ we  let $m=\min\nu$ and consider two cases.  If $\mu(h_m)\in2\bz_++1$   we use   Lemma \ref{hmueven},  followed by 
 Lemma \ref{crux}  applied to  the admissible pair $(\omega_m,\nu-\omega_m)$, to get $$(\dagger)\ \  0=h_{\omega_m,\nu-\omega_m}^{\mu}=  \sum_{\substack{\mu'\leq \nu-\omega_m}} q^{(\nu-\omega_m-\mu',\omega_m)}g_{0,\nu-\omega_m}^{\mu'} h_{\mu'+\omega_m,0}^\mu.$$
If $\mu'< \nu-\omega_m$ then $\Ht_r(\mu'+\omega_m-\mu)<\Ht_r(\nu-\mu)$ and so the inductive hypothesis  gives $$h_{\mu'+\omega_m,0}^\mu= p_{\mu'+\omega_m}^\mu= q^{(\mu'-\mu+\omega_m,\omega_m)}p_{\mu'}^{\mu-\omega_m}= q^{(\mu'-\mu+\omega_m,\omega_m)}h_{\mu',0}^{\mu-\omega_m}.$$
Substituting in $(\dagger)$, we have 
 $$0=h_{\nu,0}^\mu + q^{(\nu-\mu,\omega_m)}\sum_{\mu'<\nu-\omega_m}g_{0,\nu-\omega_m}^{\mu'} h_{\mu',0}^{\mu-\omega_m}.$$ Equation \eqref{inverse} and the inductive hypothesis applied to $\mu-\omega_m$ now give 
 $$h_{\nu,0}^\mu= q^{(\nu-\mu,\omega_m)}h_{\nu-\omega_m,0}^{\mu-\omega_m}=
 q^{(\nu-\mu,\omega_m)}p_{\nu-\omega_m}^{\mu-\omega_m}=p_\nu^\mu $$ as needed. 
 
 It remains to consider the case when  $\mu(h_m)\in 2\bz_+$ and $\nu\ne \mu$.  Since $(\nu+\omega_m)(h_m)=2$, it follows from the first part of the proof that $$(*)\ \ \ h_{\nu+\omega_m,0}^{\mu+\omega_m}= p_{\nu+\omega_m}^{\mu+\omega_m}= q^{(\nu-\mu,\omega_m)}p_\nu^\mu.$$
 Using Lemma \ref{hmueven} and Lemma \ref{crux} for the pair 
 $(\omega_m,\nu)$ gives$$(\dagger\dagger)\ \ 0=h_{\omega_m,\nu}^{\mu+\omega_m}=\sum_{\mu'\le \nu+\omega_m}q^{(\nu+\omega_m-\mu',\omega_m)} g_{0,\nu}^{\mu'-\omega_m}h_{\mu'}^{\mu+\omega_m,0}.$$ 
Since 
the inductive hypothesis applies to $h_{\mu',0}^{\mu+\omega_m}$ if $\mu'<\nu+\omega_m$
 we have $$h_{\mu',0}^{\mu+\omega_m}=p_{\mu'}^{\mu+\omega_m}= q^{(\mu'-\mu-\omega_m,\omega_m)}p_{\mu'-\omega_m}^\mu=q^{(\mu'-\mu-\omega_m,\omega_m)} h_{\mu'-\omega_m,0}^\mu.$$
 Substituting in $(\dagger\dagger)$
 and using \eqref{inverse} we get
 \begin{gather*}
0= h_{\lambda+\omega_m}^{\mu+\omega_m} + q^{(\lambda-\mu,\omega_m)}\sum_{\mu'< \lambda}g_{0,\lambda}^{\mu'} h_{\mu'}^{\mu} =h_{\lambda+\omega_m}^{\mu+\omega_m}-q^{(\lambda-\mu,\omega_m)}h_\lambda^\mu.\end{gather*} An application of equation $(*)$ completes the proof.

\end{document}